\documentclass[12pt, a4paper, twoside]{article}
\usepackage{theorem}
\ifx\pdfoutput\undefined 
\newcommand{\href}[2]{#2}
\usepackage[dvips]{graphicx}
\DeclareGraphicsExtensions{.eps}
\else 
\usepackage{color} \definecolor{webblue}{rgb}{0,0,0.6}
\usepackage[pdftex]{graphicx}
\usepackage[pdftex, colorlinks=true, linkcolor=webblue, citecolor=webblue,
pagecolor=webblue, urlcolor=webblue, pdfstartview=FitBH,
pdftitle={Homeomorphisms of the Mandelbrot Set},
pdfauthor={Wolf Jung}, bookmarksopen=true, pdfpagemode=None]{hyperref}
\DeclareGraphicsExtensions{.pdf}
\fi
\makeatletter 
\long\def\@makecaption#1#2{%
  \vskip\abovecaptionskip
  \sbox\@tempboxa{{\small\textbf{#1}: #2}}%
  \ifdim \wd\@tempboxa >\hsize
    {\small\textbf{#1}: #2}\par
  \else
    \global \@minipagefalse
    \hb@xt@\hsize{\hfil\box\@tempboxa\hfil}%
  \fi
  \vskip\belowcaptionskip}
\makeatother
\theoremstyle{break}
\theorembodyfont{\itshape}  \newtheorem{thm}{Theorem}[section]
\newtheorem{dfn}[thm]{Definition} \newtheorem{prop}[thm]{Proposition}
\theorembodyfont{\upshape}  \newtheorem{rmk}[thm]{Remark}
\newcommand{\be}{\begin{equation}} \newcommand{\ee}{\end{equation}}
\newcommand{\ban}{\begin{eqnarray}} \newcommand{\ean}{\end{eqnarray}}
\renewcommand{\hat}{\widehat} \renewcommand{\tilde}{\widetilde}
 \renewcommand{\)}{\emph{)}}
\newcommand{\eref}[1]{\emph{\ref{#1}}}
\newcommand{\mybox}{\hspace*{\fill}\rule{2mm}{2mm}}
\DeclareSymbolFont{AMSb}{U}{msb}{m}{n}
\DeclareSymbolFontAlphabet{\mathbb}{AMSb}
\def\C{\mathbb{C}}  \def\Q{\mathbb{Q}}  \def\R{\mathbb{R}}
\def\Z{\mathbb{Z}}  \def\N{\mathbb{N}}  \def\disk{\mathbb{D}}
\def\mod{\mathop{\rm mod}\nolimits}  \def\diam{\mathop{\rm diam}\nolimits}
\renewcommand{\i}{\mathrm{i}} \chardef\inodot="10
\newcommand{\eps}{\varepsilon}
\newcommand{\pre}{{\scriptscriptstyle(1)}}
\renewcommand{\r}{\mathcal{R}}  \newcommand{\E}{\mathcal{E}}
\newcommand{\G}{\mathcal{G}}  \newcommand{\K}{\mathcal{K}}
\newcommand{\M}{\mathcal{M}}  \newcommand{\sM}{{\scriptscriptstyle M}}
\renewcommand{\O}{\mathcal{O}}  \renewcommand{\P}{\mathcal{P}}
\newcommand{\bF}{\mathbf{F}}
\newcommand{\bG}{\mathbf{G}}  \newcommand{\bH}{\mathbf{H}}
\newcommand{\csd}{\hat\C\setminus\overline\disk}
\date{}
\author{Wolf Jung}
\title{Homeomorphisms of the Mandelbrot Set}
\begin{document} \maketitle
\begin{abstract}
\noindent On subsets of the Mandelbrot set, $\E_\sM\subset\M$,
homeomorphisms are constructed by quasi-conformal surgery. When the dynamics
of quadratic polynomials is changed piecewise by a combinatorial
construction, a general theorem yields the corresponding homeomorphism
$h:\E_\sM\to\E_\sM$ in the parameter plane. Each $h$ has two fixed points in
$\E_\sM\,$, and a countable family of mutually homeomorphic fundamental
domains. Possible generalizations to other families of polynomials or
rational mappings are discussed.
 \par\noindent
The homeomorphisms on subsets $\E_\sM\subset\M$ constructed by surgery are
extended to homeomorphisms of $\M$, and employed to study groups of
non-trivial homeomorphisms $h:\M\to\M$. It is shown that these groups have
the cardinality of $\R$, and they are not compact.
\end{abstract}

{\small Preprint of a paper submitted to \emph{Dynamics in the Complex
Plane}, proceedings of a symposium in honour of Bodil Branner, June 19--21
2003, Holb{\ae}k.}

\section{Introduction} \label{1}
Consider the family of complex quadratic polynomials $f_c(z):=z^2+c$. They
are parametrized by $c\in\C$, which is at the same time the critical value
of $f_c\,$, since $0$ is the critical point. The filled Julia set $\K_c$ of
$f_c$ is a compact subset of the dynamic plane. It contains all $z\in\C$
which are not attracted to $\infty$ under the iteration of $f_c\,$,
i.e.~$f_c^n(z)\not\to\infty$. The global dynamics is determined
qualitatively by the behavior of the critical point or critical value under
iteration. E.g., by a classical theorem of Fatou, $\K_c$ is connected iff
$f_c^n(c)\not\to\infty$, i.e.~$c\in\K_c\,$. The Mandelbrot set $\M$ is a
subset of the parameter plane, it contains precisely the parameters with
this property. Although it can be defined by the recursive computation of
the critical orbit, with no reference to the whole dynamic plane, most
results on $\M$ are obtained by an interplay between both planes: starting
with a subset $\E_\sM\subset\M$, employ the dynamics of $f_c$ to find a
common structure in $\K_c$ for all $c\in\E_\sM\,$. Then an analogous
structure will be found in $\E_\sM$, i.e.~in the parameter plane. This
principle has various precise formulations. Most important is its
application to external rays: these are curves in the complement of $\K_c$
(dynamic rays) or in the complement of $\M$ (parameter rays), which are
labeled by an angle $\theta\in S^1=\R/\Z$. For rational angles
$\theta\in\Q/\Z$, these rays are landing at special points in $\partial\K_c$
or $\partial\M$, respectively. See Sect.~\ref{21} for details. When rays are
landing together, the landing point is called a pinching point. It can be
used to disconnect $\K_c$ or $\M$ into well-defined components. The
structure of $\K_c\,$, as described by these pinching points, can be
understood dynamically, and then these results are transfered to the
parameter plane, to understand the structure of $\M$.

\begin{figure}[h!t!b!]
\unitlength 1mm \linethickness{0.4pt}
\begin{picture}(150, 74)
\put(10, 32){\includegraphics[width=6cm]{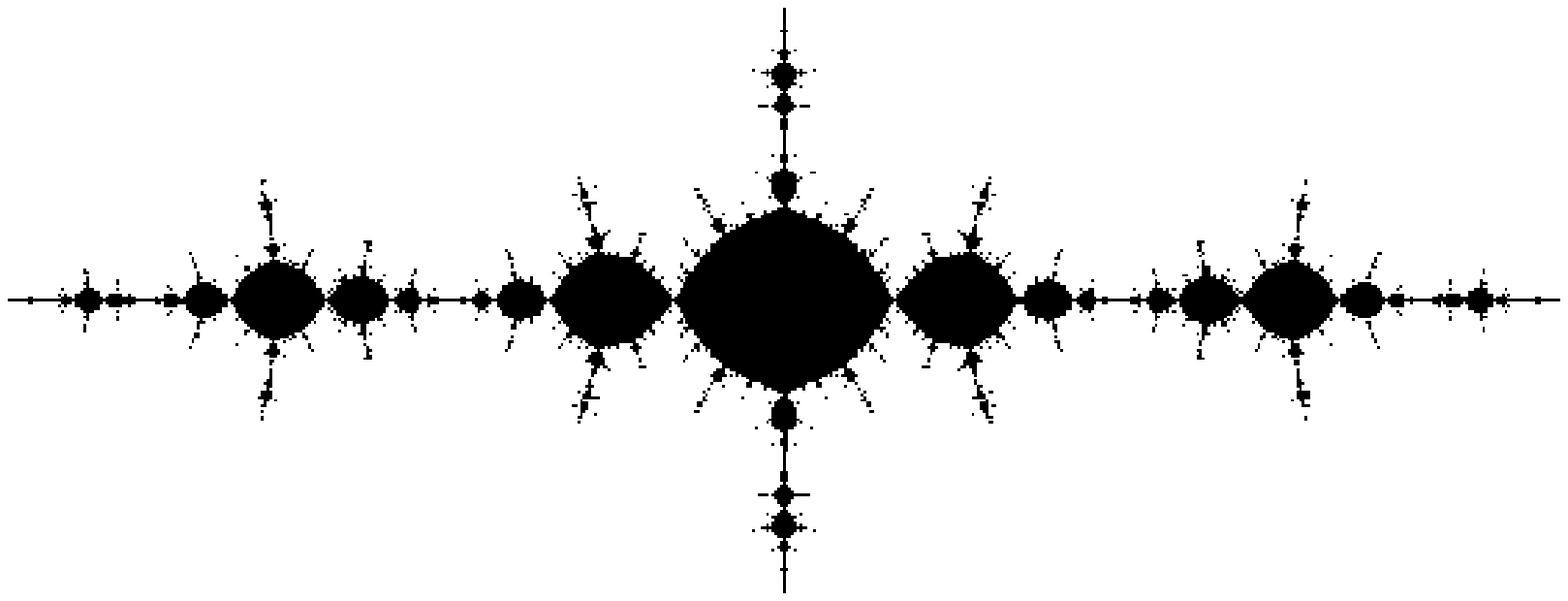}}
\put(80, 32){\includegraphics[width=6cm]{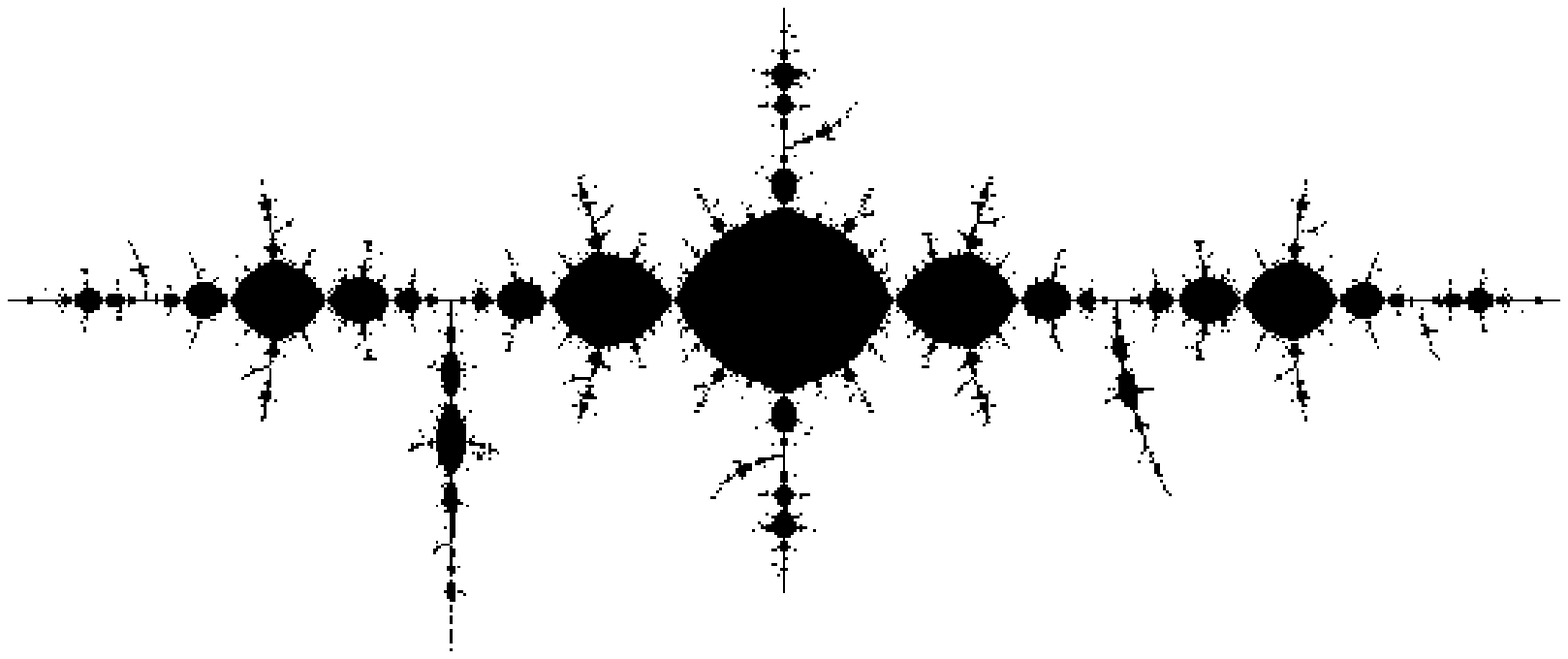}}
\put(10, 0){\includegraphics[width=6cm]{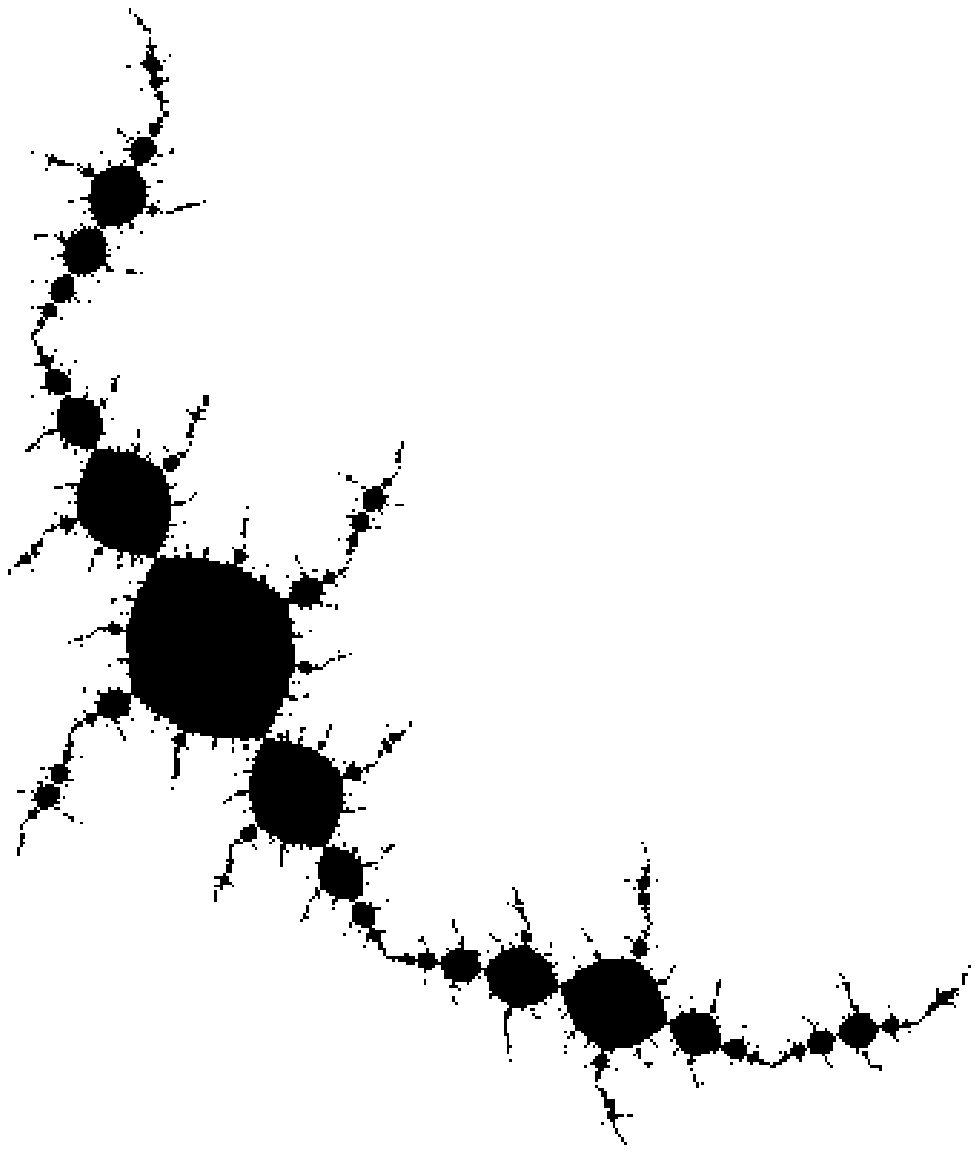}}
\put(80, 0){\includegraphics[width=6cm]{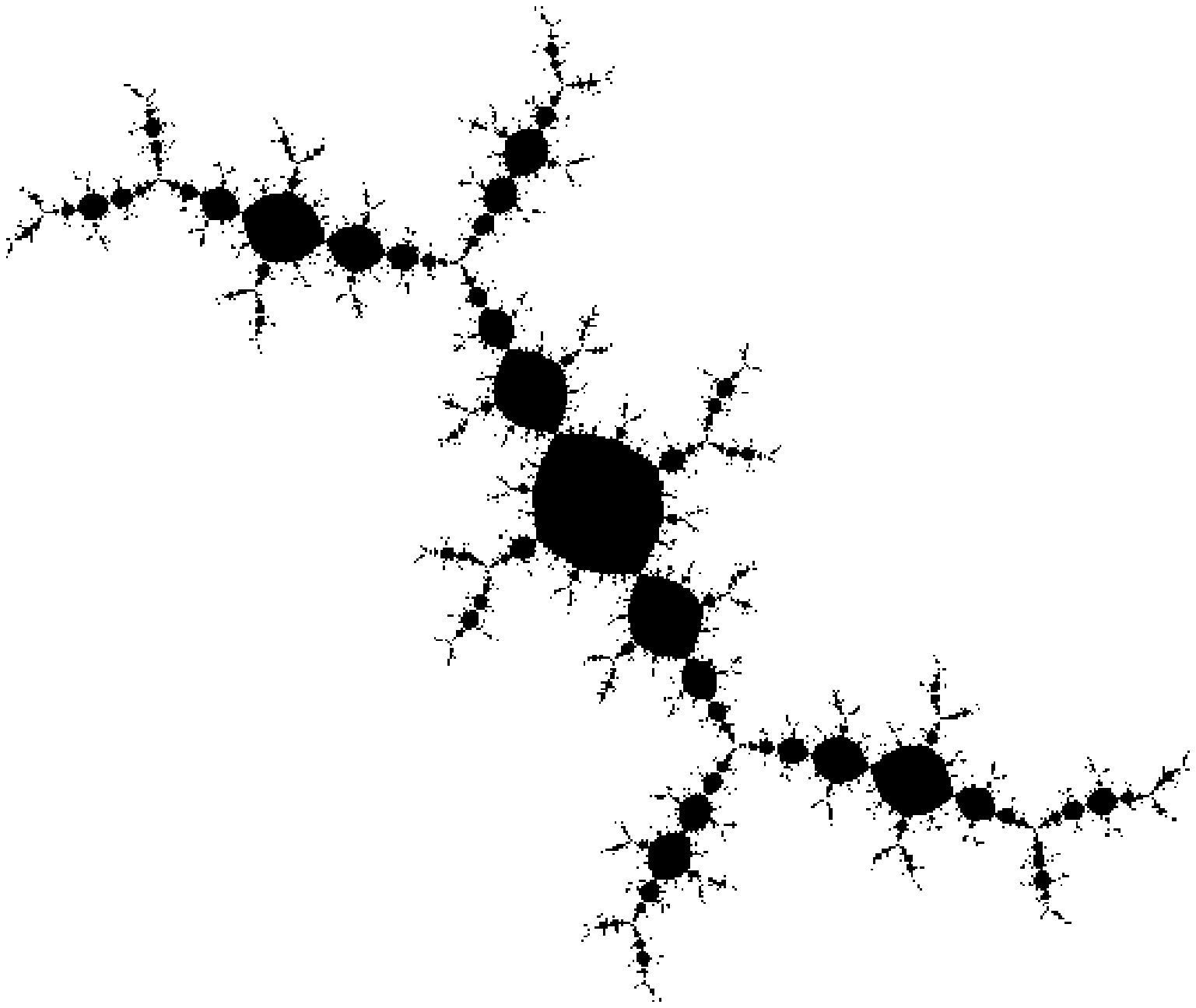}}
\thicklines
\put(110, 51){\vector(0, -1){4}}
\put(40, 47){\vector(0, 1){4}}
\put(75, 64.5){\makebox(0, 0)[cc]{$\Rightarrow$}}
\put(75, 22.5){\makebox(0, 0)[cc]{$\Leftarrow$}}
\put(113,49){\makebox(0, 0)[lc]{$\psi_c$}}
\put(37, 49){\makebox(0, 0)[rc]{$\tilde\psi_d$}}
\put(60, 56){\makebox(0, 0)[cc]{$f_c$}}
\put(90, 56){\makebox(0, 0)[cc]{$g_c$}}
\put(80, 39){\makebox(0, 0)[cc]{$f_d$}}
\put(50, 39){\makebox(0, 0)[cc]{$\tilde g_d$}}
\end{picture} \caption[Branner--Douady Surgery]{\label{Fbd}
A simulation of Branner--Douady surgery
$\Phi_A:\M_{1/2}\to\mathcal{T}\subset\M_{1/3}\,$, as explained in the text
below. In this simulation, $g_c$ and $\tilde g_d$ are defined piecewise
explicitly, and the required Riemann mappings are replaced with simple
affine mappings.}
\end{figure}

Each filled Julia set $\K_c$ is completely invariant under the corresponding
mapping $f_c\,$, and this fact explains the self-similarity of these sets.
On the other hand, when the parameter $c$ moves through the Mandelbrot set,
the corresponding Julia sets undergo an infinite number of bifurcations. By
the above principle, the local structure of $\M$ is undergoing corresponding
changes as well. But these changes may combine in such a way, that subsets
of $\M$ are mutually homeomorphic. Such homeomorphisms can be constructed by
quasi-conformal surgery. There are three basic ideas (the first and second
apply to more general situations \cite{cg}):
\begin{itemize}
\item A mapping $g$ with desired dynamics is constructed piecewise, i.e.~by
piecing together different mappings or different iterates of one mapping.
The pieces are defined e.g.~by dynamic rays landing at pinching points of
the Julia set.
\item $g$ cannot be analytic, but one constructs a quasi-conformal mapping
$\psi$ such that the composition $f=\psi\circ g\circ\psi^{-1}$ is analytic.
This is possible, when a field of infinitesimal ellipses is found that is
invariant under $g$. Then $\psi$ is constructed such that it is mapping
these ellipses to infinitesimal circles. (See Sect.~\ref{22} for the precise
definition of quasi-conformal mappings.)
\item Suppose that $f_c$ is a one-parameter family of analytic mappings,
e.g.~our quadratic polynomials, and that $g_c$ is constructed piecewise from
iterates of $f_c$ for parameters $c\in\E_\sM\subset\M$. If
 $\psi_c\circ g_c\circ\psi_c^{-1}=f_d\,$, a mapping in parameter space is
obtained from $h(c):=d$. There are techniques to show that $h$ is a
homeomorphism.
\end{itemize}
Homeomorphisms of the Mandelbrot set have been obtained in
 \cite{dhp, bd, bfl, bfe, rt, wjt}. We shall discuss the example of the
Branner--Douady homeomorphism $\Phi_A\,$, cf.~Fig.~\ref{Fbd}: parameters $c$
in the limb $\M_{1/2}$ of $\M$ are characterized by the fact, that the Julia
set $\K_c$ has two branches at the fixed point $\alpha_c$ and at its
countable family of preimages. By a piecewise construction, $f_c$ is
replaced with a new mapping $g_c\,$, such that a third branch appears at
$\alpha_c\,$, and thus at its preimages as well. This can be done by cut-
and paste techniques on a Riemann surface, or by conformal mappings between
sectors in the dynamic plane. Since $g_c$ is analytic except in some smaller
sectors, it is possible to construct an invariant ellipse field. The
corresponding quasi-conformal mapping $\psi_c$ is used to conjugate $g_c$ to
a (unique) quadratic polynomial $f_d\,$, and the mapping in parameter space
is defined by $\Phi_A(c):=d$. Now the Julia sets of $f_d$ and $g_c$ are
homeomorphic, and the dynamics are conjugate. The parameter $d$ belongs to
the limb $\M_{1/3}\,$, since the three branches of $\K_d$ at $\alpha_d$ are
permuted by $f_d$ with rotation number $1/3$. Now there is an analogous
construction of a mapping $\tilde g_d$ for
$d\in\mathcal{T}\subset\M_{1/3}\,$, which turns out to yield the inverse
mapping $\tilde\Phi_A$. The Julia set of $\tilde g_d$ has lost some arms,
and $\tilde g_d$ is conjugate to a quadratic polynomial $f_e$ again. By
showing that $f_c$ and $f_e$ are conjugate, it follows that $e=c$, thus
 $\tilde\Phi_A\circ \Phi_A=\mathrm{id}$. (The uniqueness follows from the
fact that these quasi-conformal conjugations are hybrid-equivalences,
i.e.~conformal almost everywhere on the filled Julia sets \cite{dhp}.)

\begin{figure}[h!t!b!]
\unitlength 1mm \linethickness{0.4pt}
\begin{picture}(150, 44)
\put(93, 0){\includegraphics[width=57mm]{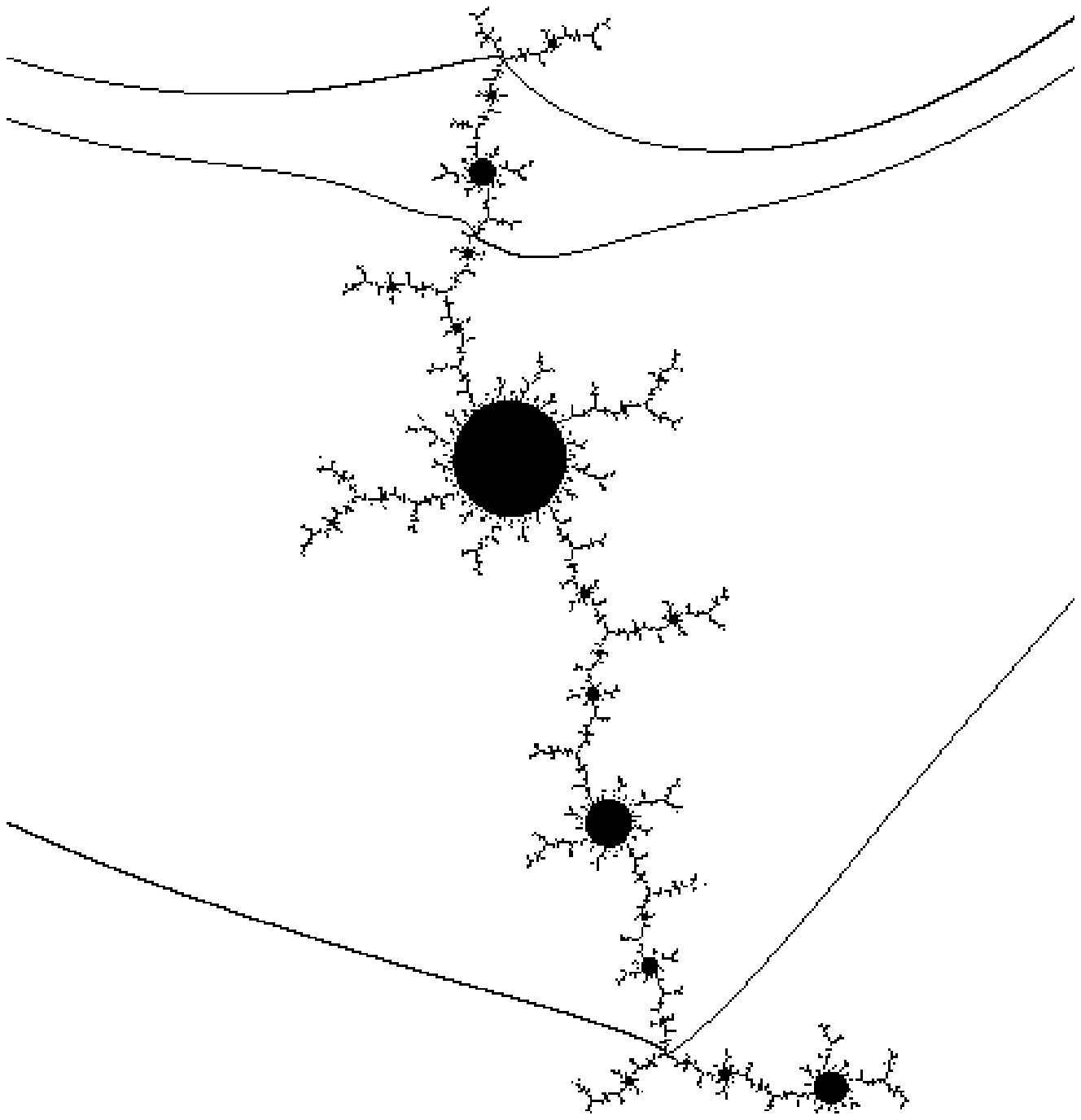}}
\put(40.5, 0){\includegraphics[width=57mm]{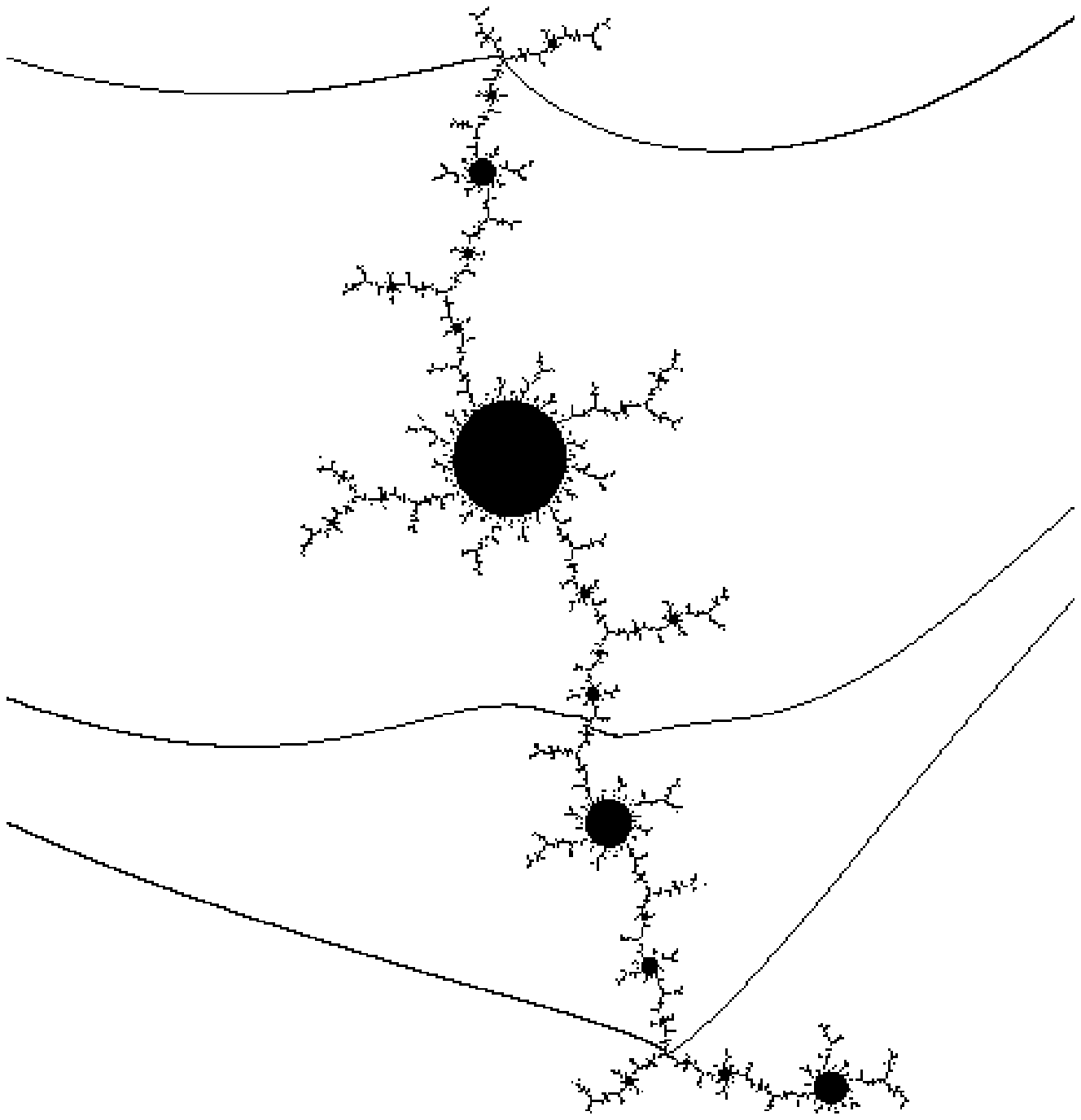}}
\put(0, 0){\includegraphics[width=57mm]{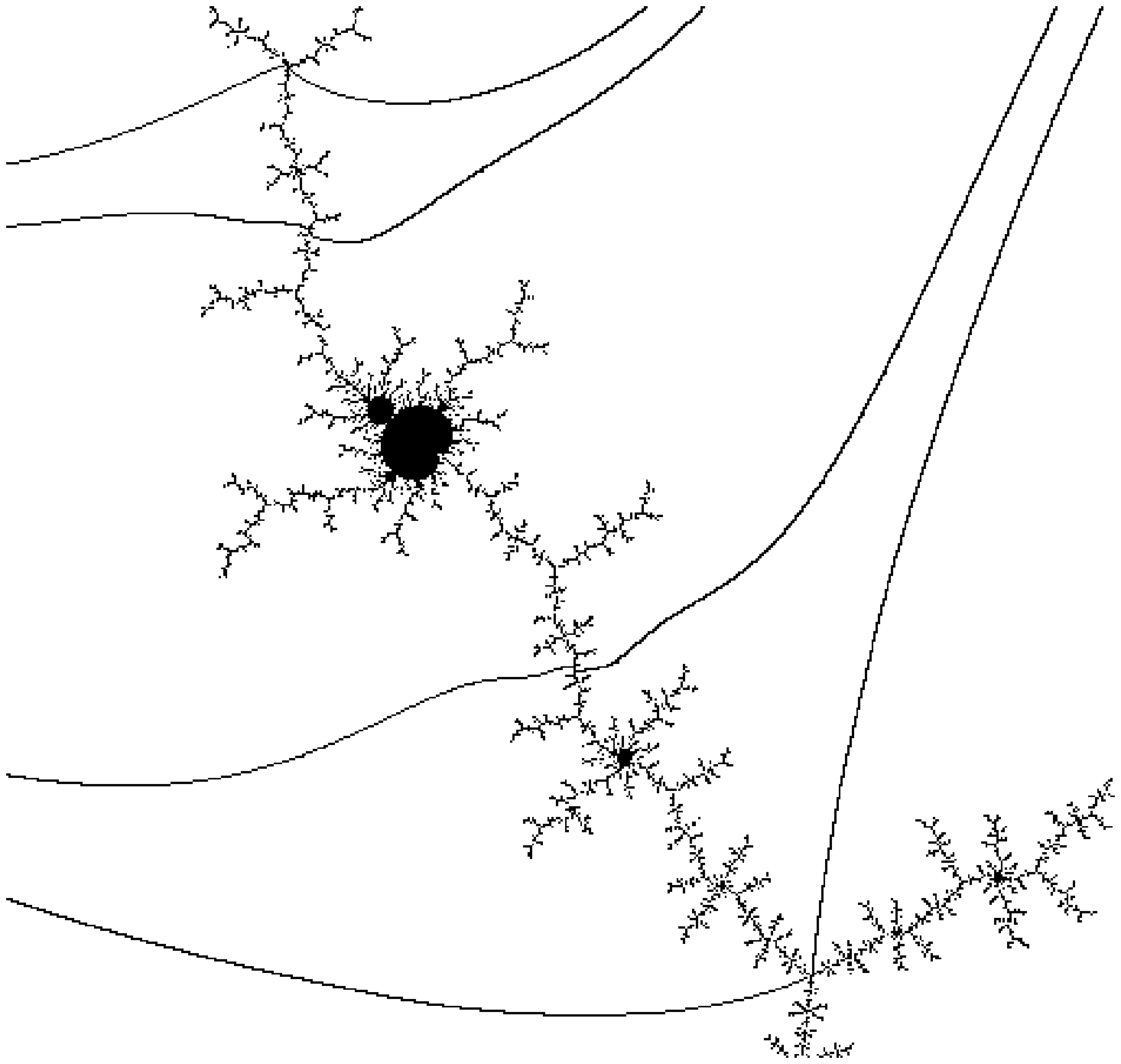}}
\put(46, 37){\makebox(0,0)[cc]{$\scriptstyle\Theta_1^-$}}
\put(36, 37){\makebox(0,0)[cc]{$\scriptstyle\Theta_2^-$}}
\put(26.5, 37){\makebox(0,0)[cc]{$\scriptstyle\Theta_3^-$}}
\put(20, 42){\makebox(0,0)[cc]{$\scriptstyle\Theta_4^-$}}
\put(3, 40.5){\makebox(0,0)[cc]{$\scriptstyle\Theta_4^+$}}
\put(3, 31.5){\makebox(0,0)[cc]{$\scriptstyle\Theta_3^+$}}
\put(3, 14){\makebox(0,0)[cc]{$\scriptstyle\Theta_2^+$}}
\put(3, 2.5){\makebox(0,0)[cc]{$\scriptstyle\Theta_1^+$}}
\put(96.5, 12){\makebox(0,0)[rc]{$\scriptstyle\Theta_1^-$}}
\put(96.5, 25){\makebox(0,0)[rc]{$\scriptstyle\Theta_2^-$}}
\put(95.5, 44){\makebox(0,0)[rt]{$\scriptstyle\Theta_4^-$}}
\put(59, 42){\makebox(0,0)[cc]{$\scriptstyle\Theta_4^+$}}
\put(59, 18){\makebox(0,0)[cc]{$\scriptstyle\Theta_2^+$}}
\put(59, 7){\makebox(0,0)[cc]{$\scriptstyle\Theta_1^+$}}
\put(149, 12){\makebox(0,0)[rc]{$\scriptstyle\Theta_1^-$}}
\put(149, 35){\makebox(0,0)[rc]{$\scriptstyle\Theta_3^-$}}
\put(148, 44){\makebox(0,0)[rt]{$\scriptstyle\Theta_4^-$}}
\put(111.5, 42){\makebox(0,0)[cc]{$\scriptstyle\Theta_4^+$}}
\put(111.5, 34){\makebox(0,0)[cc]{$\scriptstyle\Theta_3^+$}}
\put(111.5, 7){\makebox(0,0)[cc]{$\scriptstyle\Theta_1^+$}}
\put(71, 10.5){\makebox(0,0)[cc]{$V_c$}}
\put(63, 27){\makebox(0,0)[cc]{$W_c$}}
\put(114, 23){\makebox(0,0)[cc]{$\tilde V_c$}}
\put(108, 39){\makebox(0,0)[rc]{$\tilde W_c$}}
\end{picture} \caption[Definition of $g_c^\pre$]{\label{Fstrips}
Left: a parameter edge $\E_\sM\subset\M$ and the strip $\P_\sM\,$.
Middle and right: the dynamic edge $\E_c\subset\K_c$ in the strip
 $\overline{V_c\cup W_c}=\overline{\tilde V_c\cup\tilde W_c}\,$.
According to Sect.~\ref{31}, these strips are bounded by external
rays, which belong to eight angles $\Theta_i^\pm\,$. In this example,
we have $\Theta_1^-=11/56$, $\Theta_2^-=199/1008$,
$\Theta_3^-=103/504$, $\Theta_4^-=23/112$, $\Theta_4^+=29/112$,
$\Theta_3^+=131/504$, $\Theta_2^+=269/1008$, and $\Theta_1^+=15/56$.
The first-return numbers are $k_w=\tilde k_v=4,\,k_v=\tilde k_w=7$.
Now $g_c=g_c^\pre$ on $\K_c$ and $g_c^\pre=f_c\circ\eta_c\,$, with
$\eta_c=f_c^{-2}\circ(-f_c^5)=f_c^{-3}\circ(+f_c^6):V_c\to\tilde V_c\,$,
$\eta_c=f_c^{-6}\circ(-f_c^3):W_c\to\tilde W_c\,$. See also
Fig.~\ref{Fedge}.}
\end{figure}

For the homeomorphisms constructed in this paper, the mapping $g_c$ is
defined piecewise by compositions of iterates of $f_c\,$, and no cut- and
paste techniques or conformal mappings are used. Then the Julia sets of
$f_c$ and $g_c$ are the same, and no arms are lost or added in the parameter
plane either: a subset $\E_\sM\subset\M$ is defined by disconnecting $\M$ at
two pinching points, and this subset is mapped onto itself by the
homeomorphism (which is not the identity, of course). Thus a countable
family of mutually homeomorphic subsets is obtained from one construction.
General combinatorial assumptions are presented in Sect.~\ref{31}, which
allow the definition of a preliminary mapping $g_c^\pre$ analogous to the
example in Fig.~\ref{Fstrips}: it differs from $f_c$ on two strips
$V_c\,,\,W_c\,$, where it is of the form $f_c^{-n}\circ(\pm f_c^m)$.
Basically, we only need to find four strips with
 $\overline{V_c\cup W_c}=\overline{\tilde V_c\cup\tilde W_c}\,$, such that
these are mapped as $V_c\to\tilde V_c\,$, $W_c\to\tilde W_c$ by suitable
compositions of $\pm f_c^{\pm1}$.

\begin{figure}[h!t!b!]
\unitlength 1mm \linethickness{0.4pt}
\begin{picture}(150, 90)
\put(0, 45){\includegraphics[width=60mm]{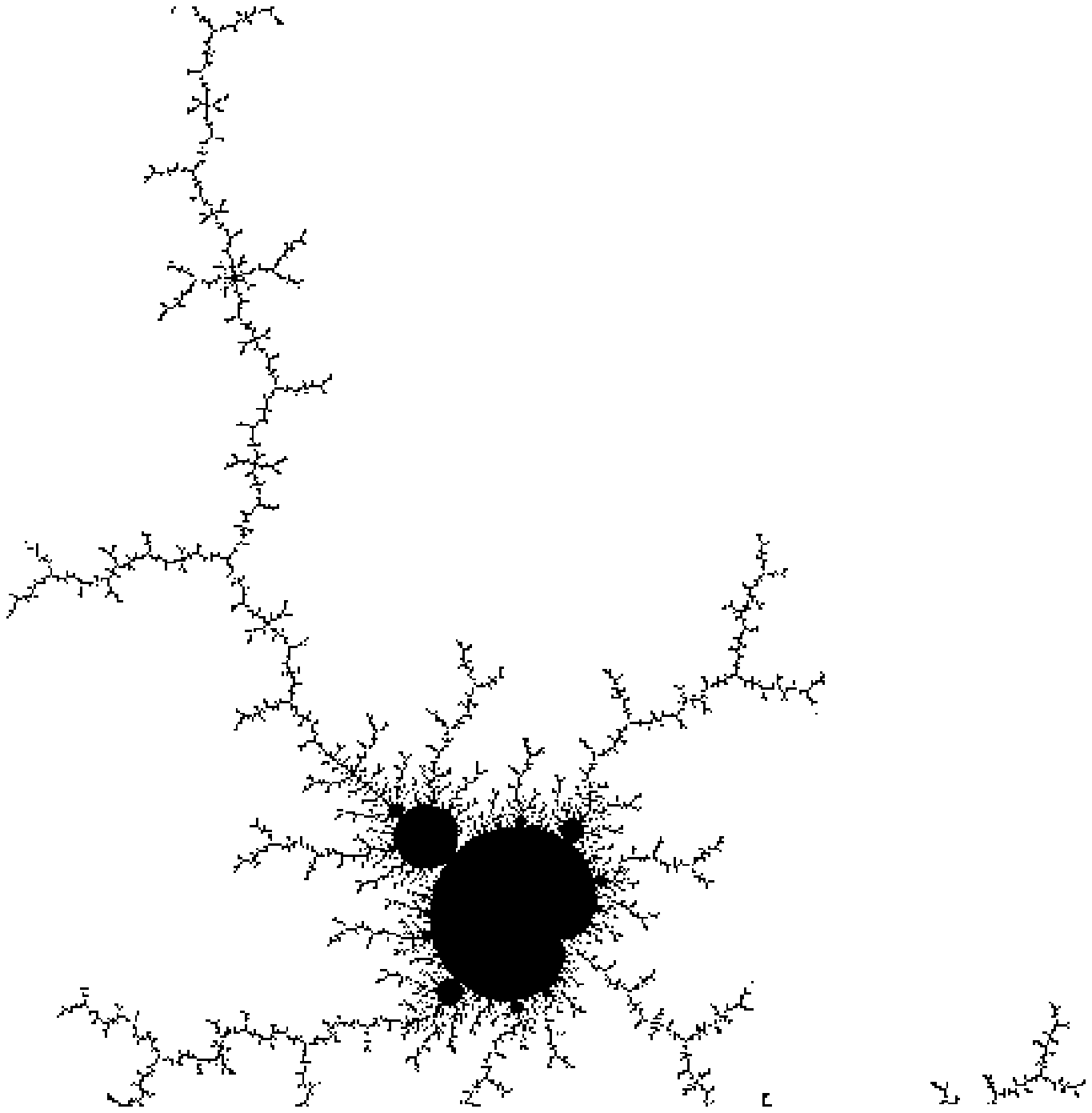}}
\put(0, 0){\includegraphics[width=60mm]{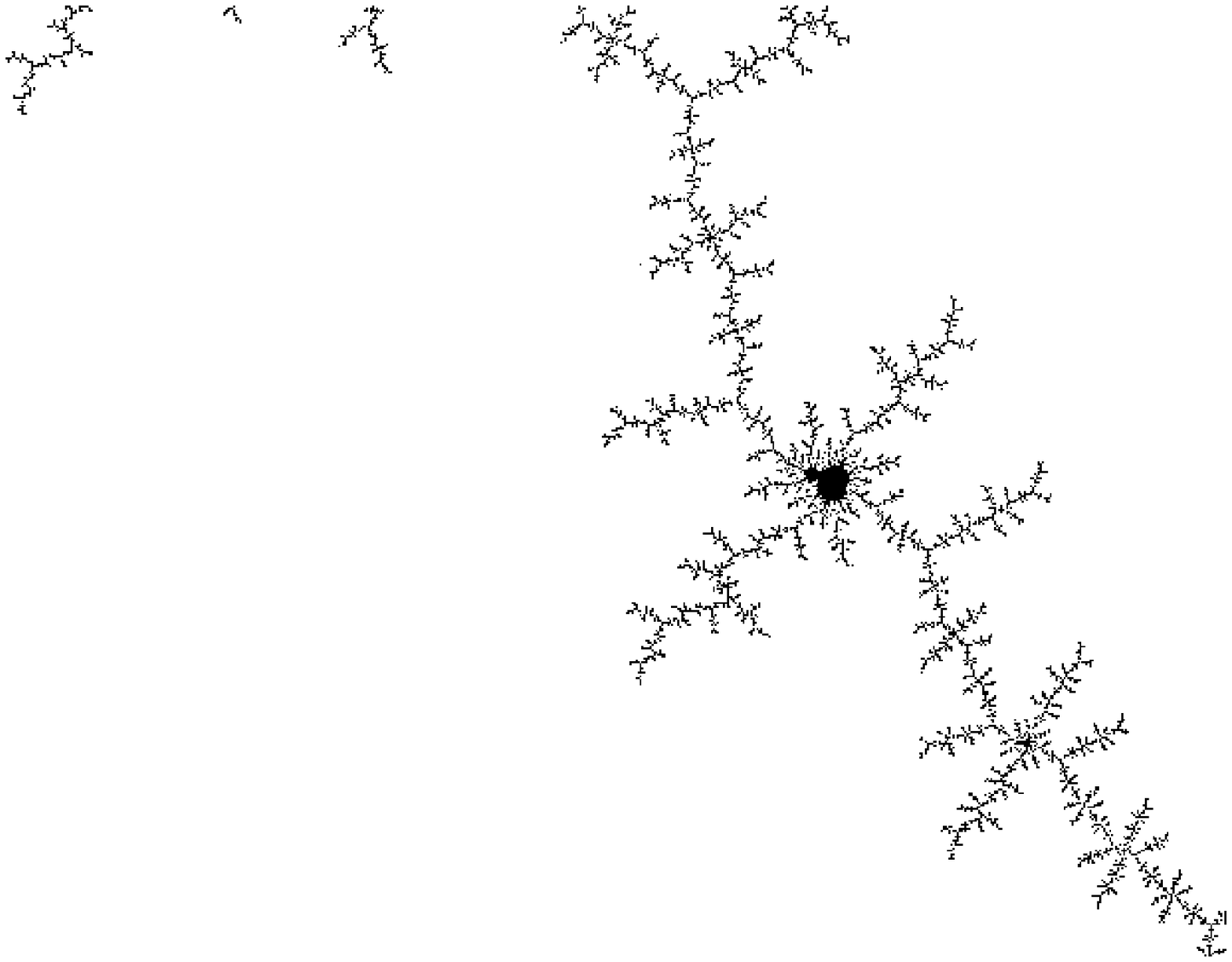}}
\put(110, 0){\includegraphics[width=40mm]{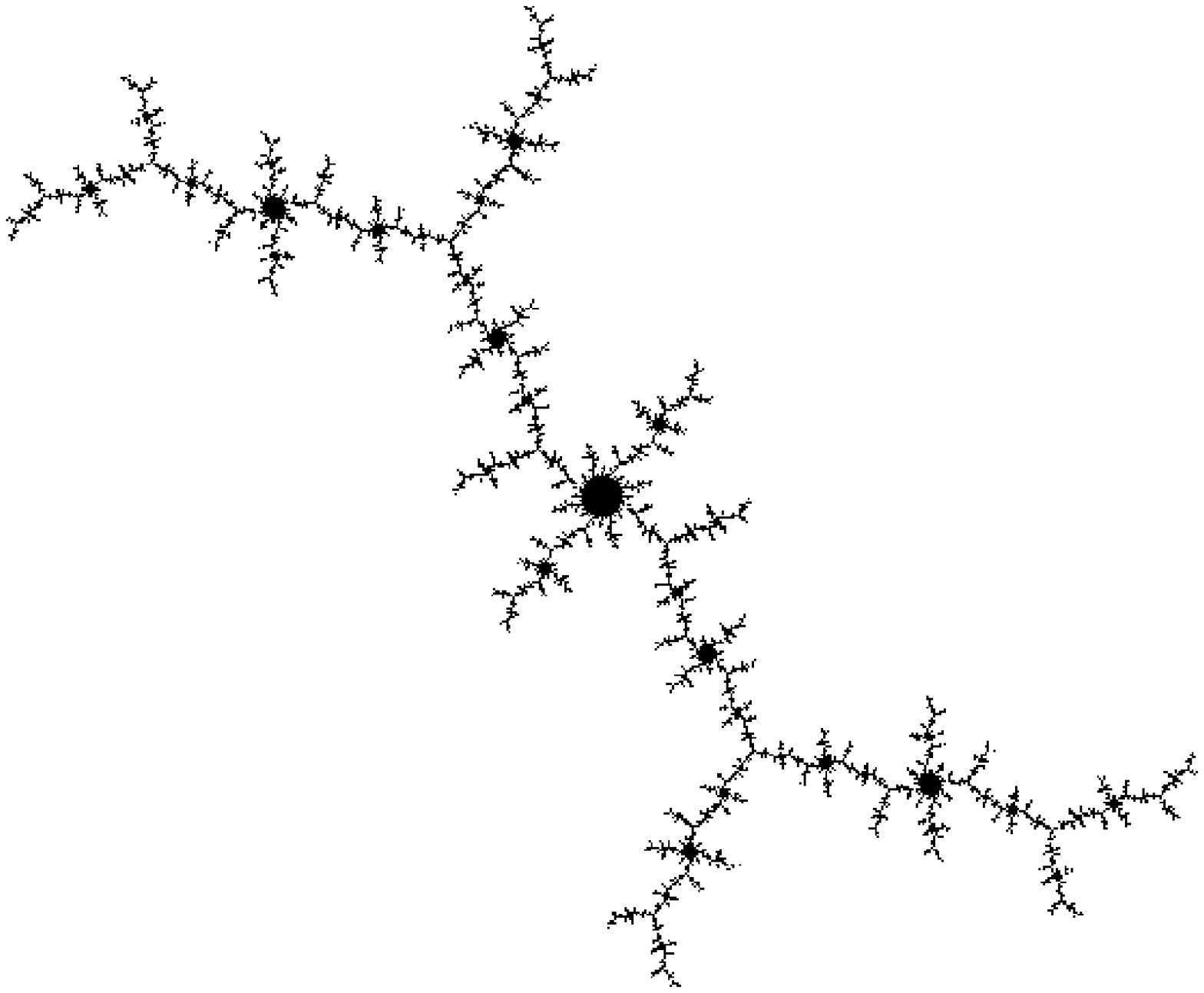}}
\put(90, 30){\includegraphics[width=40mm]{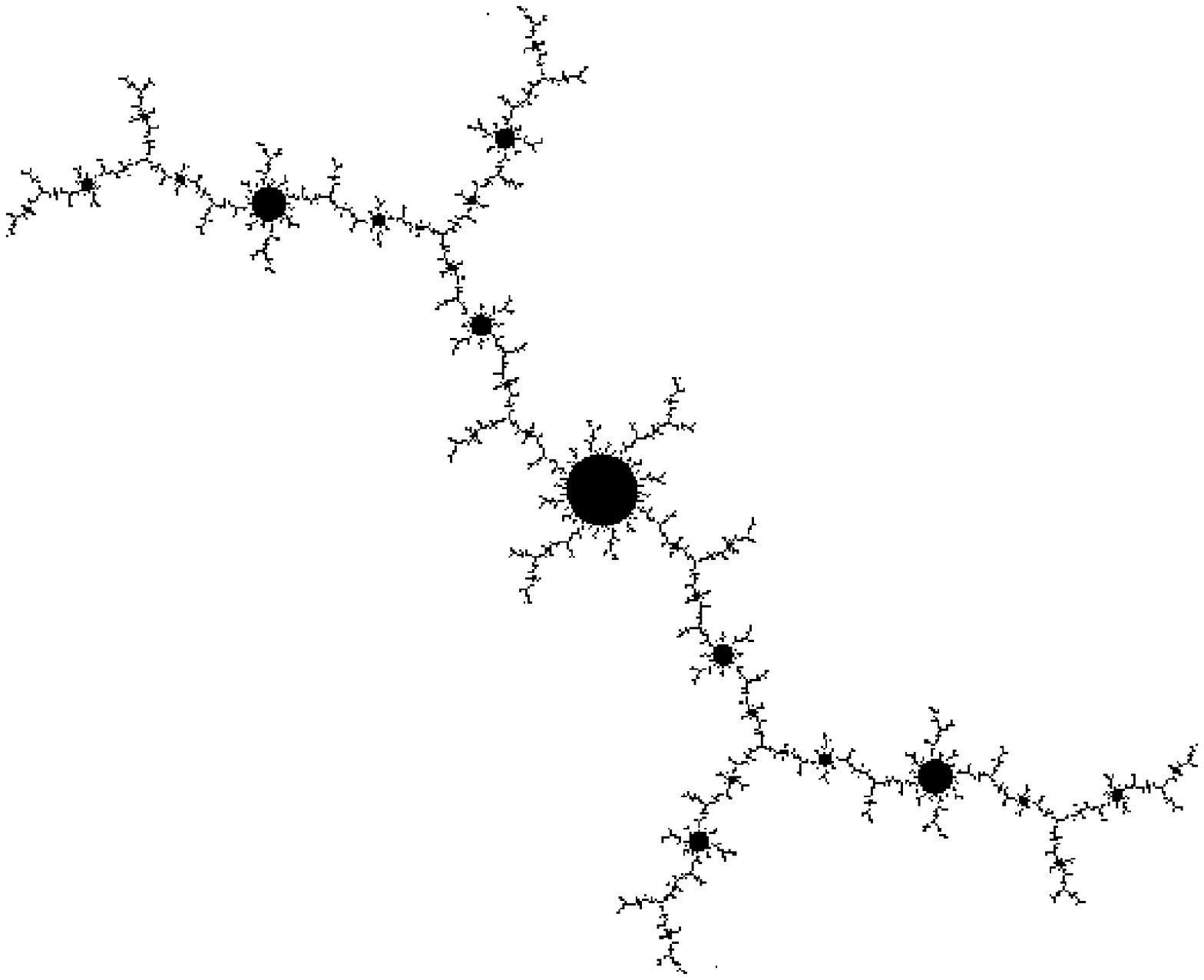}}
\put(70, 60){\includegraphics[width=40mm]{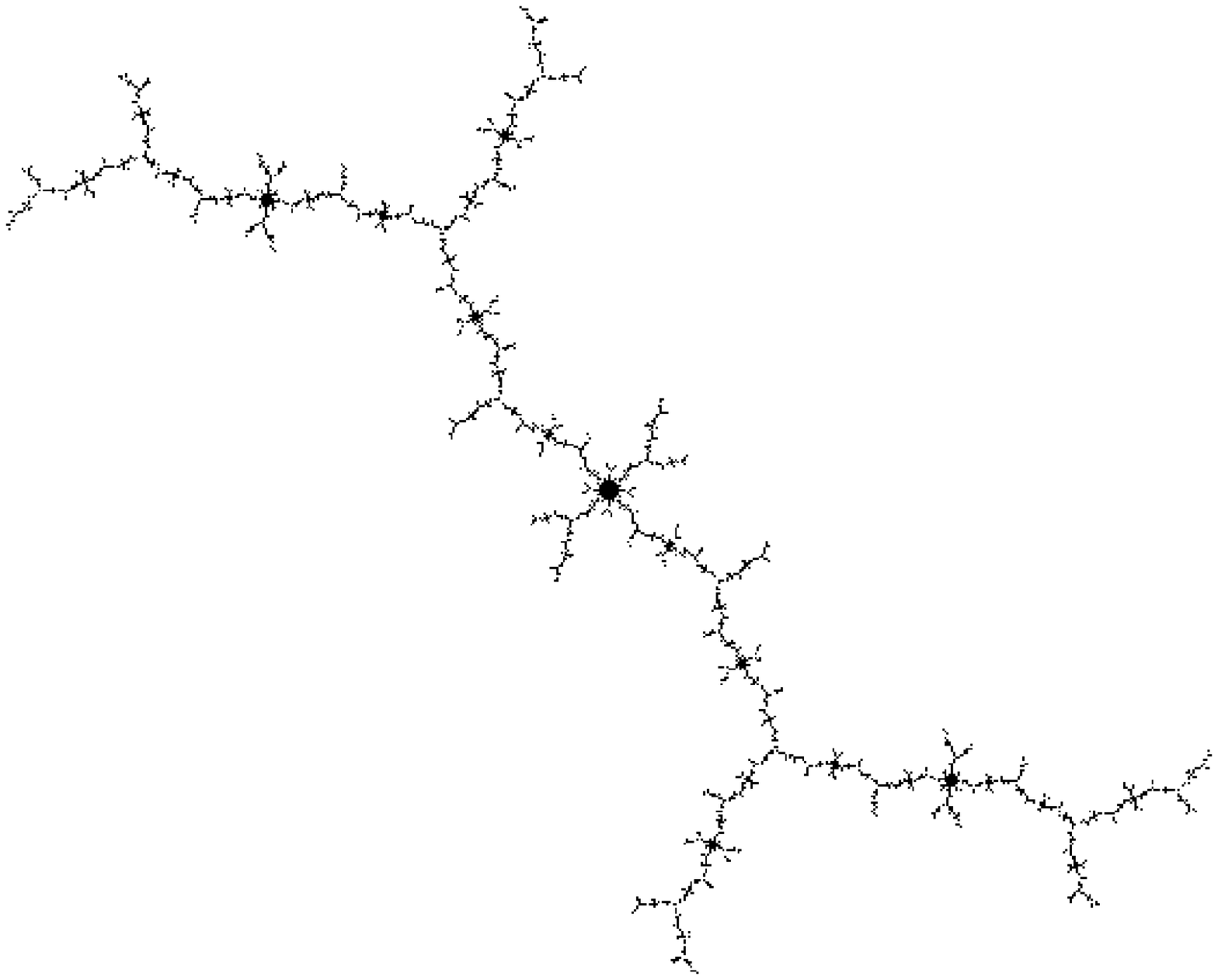}}
\thicklines
\put(65, 5){\makebox(0,0)[cc]{$a$}}
\put(6, 89){\makebox(0,0)[ct]{$b$}}
\put(63, 34){\makebox(0,0)[cc]{$c_7'$}}
\put(47, 58){\makebox(0,0)[cc]{$c_4$}}
\put(31, 82){\makebox(0,0)[cc]{$c_7''$}}
\put(60, 38.5){\vector(-2, 3){10}}
\put(44, 61.5){\vector(-2, 3){10}}
\put(60, 48){\makebox(0,0)[cc]{$h$}}
\put(44, 72){\makebox(0,0)[cc]{$h$}}
\put(119, 10){\makebox(0,0)[cc]{$\K_{c_7'}$}}
\put(99, 40){\makebox(0,0)[cc]{$\K_{c_4}$}}
\put(79, 70){\makebox(0,0)[cc]{$\K_{c_7''}$}}
\put(138, 29){\vector(-2, 3){12}}
\put(118, 59){\vector(-2, 3){12}}
\put(137, 40){\makebox(0,0)[cc]{$\psi_{c_7'}$}}
\put(117, 70){\makebox(0,0)[cc]{$\psi_{c_4}$}}
\end{picture} \caption[Homeomorphism on an Edge]{\label{Fedge}
Left: the parameter edge $\E_\sM$ from $a:=\gamma_\sM(11/56)$ to
$b:=\gamma_\sM(23/112)$, the same as in Fig.~\ref{Fstrips}. The
homeomorphism $h:\E_\sM\to\E_\sM$ is expanding at $a$ and contracting at
$b$. The centers of periods $4$ and $7$ are mapped as
 $h:c_7'\mapsto c_4\mapsto c_7''\,$. Right: the filled Julia sets for
$c_7'\,$, $c_4\,$, and $c_7''$ are quasi-conformally homeomorphic.
($\E_c\subset\K_c$ is barely visible in the top right corner.)}
\end{figure}

\begin{thm}[Construction and Properties of $h$]\label{Th}
$1.$ Given the combinatorial construction of $\E_\sM\subset\M$ and
$g_c^\pre$ for $c\in\E_\sM$ according to Def.~\eref{Dg}, there is a family
of ``quasi-quadratic'' mappings $g_c$ coinciding with $g_c^\pre$ on the
filled Julia sets $\K_c\,$. These are hybrid-equivalent to unique quadratic
polynomials.

$2.$ The mapping $h:\E_\sM\to\E_\sM$ in parameter space is defined as
follows: for $c\in\E_\sM\,$, find the polynomial $f_d$ that is
hybrid-equivalent to $g_c\,$, and set $h(c):=d$. It does not depend on the
precise choice of $g_c$ \(only on the combinatorial definition of
$g_c^\pre$\). Now $h$ is a homeomorphism, and analytic in the interior of
$\E_\sM\,$.

$3.$ $h$ is a non-trivial homeomorphism of $\E_\sM$ onto itself, fixing the
vertices $a$ and $b$. $h$ and $h^{-1}$ are H\"older continuous at
Misiurewicz points and Lipschitz continuous at $a$ and $b$. Moreover, $h$ is
expanding at $a$ and contracting at $b$, cf.~Fig.~\eref{Fedge}: for
$c\in\E_\sM\setminus\{a,\,b\}$ we have $h^n(c)\to b$ as $n\to\infty$, and
$h^{-n}(c)\to a$. There is a countable family of mutually homeomorphic
fundamental domains.

$4.$ $h$ extends to a homeomorphism between strips,
$h:\P_\sM\to\tilde\P_\sM\,$, which is quasi-conformal in the exterior of
$\M$.
\end{thm}

The power of Thm.~\ref{Th} lies in turning combinatorial data into
homeomorphisms. The creative step remaining is to find eight angles
$\Theta_i^\pm$, such that there are compositions of $\pm f_c^{\pm1}$ mapping
$V_c\to\tilde V_c$ and $W_c\to\tilde W_c\,$. When this is done, the
existence of a corresponding homeomorphisms is guaranteed by the theorem.

In Sections~\ref{2} and~\ref{3}, basic properties of $\M$ and of
quasi-conformal mappings are recalled, and the proof of Thm.~\ref{Th} is
sketched by constructing $g_c$ and $h$. Related results from the author's
thesis \cite{wjt} are summarized in Sect.~\ref{4}. These include more
examples of homeomorphisms, constructed at chosen Misiurewicz points or on
``edges,'' and the combinatorial description of homeomorphisms. When a
one-parameter family of polynomials is defined by critical relations,
homeomorphisms in parameter space can be obtained by analogous techniques.
The same applies e.g.~to the rational mappings arising in Newton's method
for cubic polynomials.

H.~Kriete has suggested that the homeomorphisms $h:\E_\sM\to\E_\sM$
constructed by this kind of surgery extend to homeomorphisms of $\M$. Thus
they can be used to study the homeomorphism group of $\M$, answering a
question by K.~Keller. In Sect.~\ref{5}, some possible definitions for
groups of non-trivial homeomorphisms are discussed, and some properties are
obtained by combining two tools: the characterization of homeomorphisms by
permutations of hyperbolic components, and the composition of homeomorphisms
constructed by surgery. The groups are not compact, and the groups of
non-trivial homeomorphisms have the cardinality of $\R$.

\subsection*{Acknowledgment}
Many people have contributed to this work by inspiring discussions and
helpful suggestions. I wish to thank in particular Mohamed Barakat, Walter
Bergweiler, Volker Enss, N\'uria Fagella, John Hamal Hubbard, Gerhard Jank,
Karsten Keller, Hartje Kriete, Fernando Lled\'o, Olaf Post, Johannes Riedl,
and Dierk Schleicher.

I am especially happy to contribute this paper to the proceedings of a
conference in honor of Bodil Branner, since I have learned surgery from her
papers \cite{bd, bfl}.

\newpage
\section{Background} \label{2}
Our main tools are the landing properties of external rays, and sending an
ellipse field to circles by a quasi-conformal mapping.

\subsection{The Mandelbrot Set} \label{21}
$f_c(z)=z^2+c$ has a superattracting fixed point at
$\infty\in\hat\C:=\C\cup\{\infty\}$. The unique Boettcher conjugation is
conjugating $f_c$ to $F(z):=z^2$, $\Phi_c\circ f_c=F\circ\Phi_c$ in a
neighborhood of $\infty$. If the critical point $0$, or the critical value
$c$, does not escape to $\infty$, then $\K_c$ is connected \cite{cg}, and
the parameter $c$ belongs to the Mandelbrot set $\M$ by definition. Then
$\Phi_c$ extends to a conformal mapping $\Phi_c:\hat\C\setminus\K_c\to\csd$,
where $\overline\disk$ is the closed unit disk. Dynamic rays $\r_c(\theta)$
are defined as preimages of straight rays
$\r(\theta)=\{z\,|\,1<|z|<\infty,\,\arg(z)=2\pi\theta\}$ under $\Phi_c\,$.
Theses curves in the complement of $\K_c$ may land at a point in
$\partial\K_c\,$, or accumulate at the boundary without landing. But they
are always landing when $\theta$ is rational. Each rational angle $\theta$
is periodic or preperiodic under doubling ($\mod1$). In the former case, the
dynamic ray $\r_c(\theta)$ is landing at a periodic point
$z=\gamma_c(\theta)\in\partial\K_c\,$, and at a preperiodic point in the
latter case. These properties are understood from
$f_c(\r_c(\theta))=\r_c(2\theta)$, since $\arg(F(z))=2\arg(z)$.

The Mandelbrot set is compact, connected, and full, and the conformal
mapping $\Phi_\sM:\hat\C\setminus\M\to\csd$ is given by
$\Phi_\sM(c):=\Phi_c(c)$. (When $c\notin\M$, $\K_c$ is totally disconnected
and $\Phi_c$ is not defined in all of its complement, but it is well-defined
at the critical value.) Parameter rays $\r_\sM(\theta)$ are defined as
preimages of straight rays under $\Phi_\sM\,$. Their landing properties are
obtained e.g.~in \cite{ser}: each rational ray $\r_\sM(\theta)$ is landing
at a point $c=\gamma_\sM(\theta)\in\partial\M$. When $\theta$ is
preperiodic, then the critical value $c$ of $f_c$ is preperiodic, and the
parameter $c$ is called a Misiurewicz point. The critical value $c\in\K_c$
has the same external angles as the parameter $c\in\M$. When $\theta$ is
periodic, then $c$ is the root of a hyperbolic component (see below). Both
in the dynamic plane of $f_c$ and in the parameter plane, the landing points
of two or more rational rays are called pinching points. They are used to
disconnect these sets into well-defined components, which are described
combinatorially by rational numbers. Their structure is obtained from the
dynamics, and transfered to the parameter plane. Pinching points with more
than two branches are branch points.

Hyperbolic components of $\M$ consist of parameters, such that
the corresponding polynomial has an attracting cycle. The root is the
parameter on the boundary, such that the cycle has multiplier $1$. The
boundary of a hyperbolic component contains a dense set of roots of
satellite components. Each hyperbolic component has a unique center, where
the corresponding cycle is superattracting. Centers or roots are dense at/in
$\partial\M$.

\newpage
\subsection{Quasi-Conformal Mappings} \label{22}
An orientation-preserving homeomorphism $\psi$ between domains in $\hat\C$
is \emph{$K$-quasi-conformal}, $1\le K<\infty$, if it has two properties:
\begin{itemize}
 \item It is weakly differentiable, so that its differential
$d\psi=\partial\psi\,dz+\overline\partial\psi\,d\overline z$ is defined
almost everywhere. This linear map is sending certain ellipses in the
tangent space to circles. The Beltrami coefficient
$\mu:=\overline\partial\psi/\partial\psi$ is defined almost everywhere. It
encodes the direction and the dilatation ratio of the semi-axes for the
ellipse field \cite{cg}.
 \item The dilatation ratio is bounded globally by $K$, or
 $|\mu(z)|\le(K-1)/(K+1)$ almost everywhere.
\end{itemize}
The chain rule for derivatives is satisfied for the composition of
quasi-conformal mappings, and a $1$-quasi-conformal mapping is conformal.
Quasi-conformal mappings are absolutely continuous, H\"older continuous, and
have nice properties regarding e.g.~boundary behavior or normal families
\cite{lvqk}. Given a measurable ellipse field (Beltrami coefficient) $\mu$
with $|\mu(z)|\le m<1$ almost everywhere, the Beltrami differential equation
$\overline\partial\psi=\mu\partial\psi$ on $\hat\C$ has a unique solution
with the normalization $\psi(z)=z+o(1)$ as $z\to\infty$. The dependence on
parameters is described by the Ahlfors--Bers Theorem \cite{cg}, which is
behind some of our arguments, but will not be used explicitly here.

A \emph{$K$-quasi-regular} mapping is locally $K$-quasi-conformal except for
critical points, but it need not be injective globally. In Sect.~\ref{33},
we will have a quasi-regular mapping $g$, such that all iterates are
$K$-quasi-regular, and analytic in a neighborhood of $\infty$. Then a
$g$-invariant field of infinitesimal ellipses is obtained as follows: it
consists of circles in a neighborhood of $\infty$, i.e.~$\mu(z)=0$ there,
and it is pulled back with iterates of $g$. Now $\psi$ shall solve the
corresponding Beltrami equation, i.e.~send these ellipses to circles. By the
chain rule, $f:=\psi\circ g\circ\psi^{-1}$ is mapping almost every
infinitesimal circle to a circle, thus it is analytic.

\section{Quasi-Conformal Surgery} \label{3}
As soon as the combinatorial assumptions on $g_c^\pre$ given here are
satisfied, Thm.~\ref{Th} yields a corresponding homeomorphism of $\E_\sM\,$.
After formulating these general assumptions, the proof is sketched by
constructing the quasi-quadratic mapping $g_c$ and the homeomorphism $h$.
For some details, the reader will be referred to \cite{wjt}.

\subsection{Combinatorial Setting} \label{31}
The following definitions may be illustrated by the example in
Fig.~\ref{Fstrips}. Further examples are mentioned in
Sects.~\ref{42}--\ref{44}. When four parameter rays are landing in pairs at
two pinching points of $\M$, this defines a strip in the parameter plane.
Analogously, four dynamic rays define a strip in the dynamic plane. Our
assumptions are formulated in terms of eight preperiodic angles
\be 0<\Theta_1^-<\Theta_2^-<\Theta_3^-<\Theta_4^-<
 \Theta_4^+<\Theta_3^+<\Theta_2^+<\Theta_1^+<1 \ . \ee
\begin{itemize}
 \item The Misiurewicz points
 $a:=\gamma_\sM(\Theta_1^-)=\gamma_\sM(\Theta_1^+)\neq
 \gamma_\sM(\Theta_4^-)=\gamma_\sM(\Theta_4^+)=:b$ mark a compact,
connected, full subset $\E_\sM\subset\M$: $\E_\sM=\P_\sM\cap\M$, where
$\P_\sM$ is the closed strip bounded by the four parameter rays
$\r_\sM(\Theta_1^\pm),\,\r_\sM(\Theta_4^\pm)$.
 \item For all $c\in\E_\sM$, the eight dynamic rays $\r_c(\Theta_i^\pm)$
shall be landing in pairs at four distinct points,
i.e.~$\gamma_c(\Theta_i^-)=\gamma_c(\Theta_i^+)$. (Equivalently, they are
landing in this pattern for one $c_0\in\E_\sM\,$, and none of the eight
angles is returning to $(\Theta_1^-,\,\Theta_1^+)$ under doubling $\mod1$.)
Four open strips are defined as follows, cf.~Fig.~\ref{Fstrips}:
 $V_c$ is bounded by $\r_c(\Theta_1^\pm)$ and $\r_c(\Theta_2^\pm)$,
 $W_c$ is bounded by $\r_c(\Theta_2^\pm)$ and $\r_c(\Theta_4^\pm)$,
 $\tilde V_c$ is bounded by $\r_c(\Theta_1^\pm)$ and $\r_c(\Theta_3^\pm)$,
 $\tilde W_c$ is bounded by $\r_c(\Theta_3^\pm)$ and $\r_c(\Theta_4^\pm)$.
$\E_c\subset\K_c$ is defined as the intersection of $\K_c$ with the closed
strip $\overline{V_c\cup W_c}=\overline{\tilde V_c\cup\tilde W_c}\,$. Thus
for parameters $c\in\E_\sM\,$, the critical value $c$ satisfies
$c\in\E_c\,$.
 \item The first-return number $k_v$ is the smallest integer $k>0$, such
that $f_c^k(V_c)$ meets (covers) $\E_c\,$. Equivalently, it is the largest
integer $k>0$, such that $f_c^{k-1}$ is injective on $V_c\,$. Define
 $k_w\,,\,\tilde k_v\,,\,\tilde k_w$ analogously. They are independent of
$c\in\E_\sM\,$. Now the main assumption on the dynamics, which makes finding
the angles non-trivial, is that there is a (fixed) choice of signs with
 $f_c^{k_v-1}(V_c)=\pm f_c^{\tilde k_v-1}(\tilde V_c)$ and
 $f_c^{k_w-1}(W_c)=\pm f_c^{\tilde k_w-1}(\tilde W_c)$. The ``orientation''
is respected, i.e.~with $z_i:=\gamma_c(\Theta_i^\pm)$ we have
 e.g.~$f_c^{k_v-1}(z_1)=\pm f_c^{\tilde k_v-1}(z_1)$ and
 $f_c^{k_v-1}(z_2)=\pm f_c^{\tilde k_v-1}(z_3)$.
\end{itemize}

If $a$ or $b$ is a branch point of $\M$, the last assumption implies that
$\E_\sM$ is contained in a single branch, i.e.~$\E_\sM\setminus\{a,\,b\}$ is
a connected component of $\M\setminus\{a,\,b\}$.

\begin{dfn}[Preliminary Mapping $g_c^\pre$] \label{Dg}
Under these assumptions, with the unique choices of signs in the two strips,
define
 $\eta_c:=f_c^{-(\tilde k_v-1)}\circ(\pm f_c^{k_v-1}):V_c\to \tilde V_c\,$,
 $\eta_c:=f_c^{-(\tilde k_w-1)}\circ(\pm f_c^{k_w-1}):W_c\to \tilde W_c\,$,
and $\eta_c:=\mathrm{id}$ on $\C\setminus\overline{V_c\cup W_c}$ for
$c\in\E_\sM$. Then define $g_c^\pre:=f_c\circ\eta_c$ and
 $\tilde g_c^\pre:=f_c\circ\eta_c^{-1}$.
\end{dfn}

The three mappings are holomorphic and defined piecewise, thus they cannot
be extended continuously. Each has ``shift discontinuities'' on six dynamic
rays: e.g., consider $z_0\in\r_c(\Theta_2^-)$, $(z_n')\subset V_c$ and
 $(z_n'')\subset W_c$ with $z_n'\to z_0$ and $z_n''\to z_0\,$, then
 $\lim g_c^\pre(z_n')$ and $\lim g_c^\pre(z_n'')$ both exist and belong to
$\r_c(\Theta_3^-)$, but they are shifted relative to each other along this
ray. Neglecting these rays, $g_c^\pre$ and $\tilde g_c^\pre$ are proper of
degree $2$. In the following section, $g_c^\pre$ will be replaced with a
smooth mapping $g_c\,$, which is used to construct the homeomorphism $h$.
Analogously, $\tilde g_c^\pre$ yields $\tilde h=h^{-1}$.

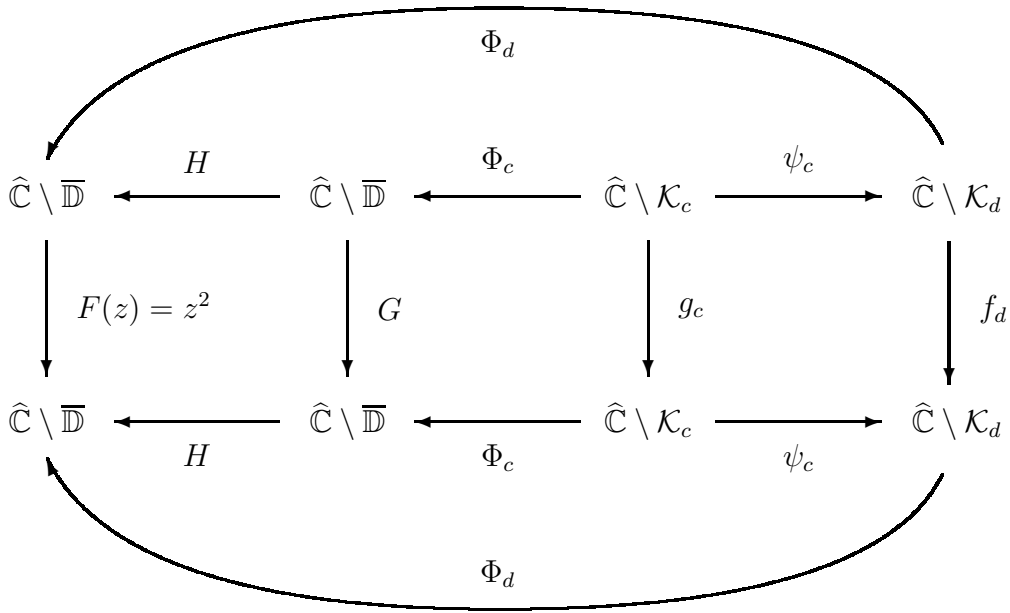
\begin{figure}[h!t!b!]
\unitlength 1mm \linethickness{0.4pt}
\begin{picture}(150, 80)\put(10, 0){\begin{picture}(130, 80)
\thicklines
\put(5, 55){\makebox(0,0)[cc]{$\csd$}}
\put(5, 25){\makebox(0,0)[cc]{$\csd$}}
\put(9, 40){\makebox(0,0)[lc]{$F(z)=z^2$}} \put(5, 49){\vector(0, -1){18}}
\put(45, 55){\makebox(0,0)[cc]{$\csd$}}
\put(45, 25){\makebox(0,0)[cc]{$\csd$}}
\put(49, 40){\makebox(0,0)[lc]{$G$}} \put(45, 49){\vector(0, -1){18}}
\put(85, 55){\makebox(0,0)[cc]{$\hat\C\setminus\K_c$}}
\put(85, 25){\makebox(0,0)[cc]{$\hat\C\setminus\K_c$}}
\put(89, 40){\makebox(0,0)[lc]{$g_c$}} \put(85, 49){\vector(0, -1){18}}
\put(126, 55){\makebox(0,0)[cc]{$\hat\C\setminus\K_d$}}
\put(126, 25){\makebox(0,0)[cc]{$\hat\C\setminus\K_d$}}
\put(129, 40){\makebox(0,0)[lc]{$f_d$}} \put(125, 49){\vector(0, -1){19}}
\put(25, 58){\makebox(0,0)[cb]{$H$}} \put(36, 55){\vector(-1, 0){22}}
\put(25, 22){\makebox(0,0)[ct]{$H$}} \put(36, 25){\vector(-1, 0){22}}
\put(65, 58){\makebox(0,0)[cb]{$\Phi_c$}} \put(76, 55){\vector(-1, 0){22}}
\put(65, 22){\makebox(0,0)[ct]{$\Phi_c$}} \put(76, 25){\vector(-1, 0){22}}
\put(105, 58){\makebox(0,0)[cb]{$\psi_c$}} \put(94, 55){\vector(1, 0){22}}
\put(105, 22){\makebox(0,0)[ct]{$\psi_c$}} \put(94, 25){\vector(1, 0){22}}
\qbezier(5, 60)(15, 80)(65, 80) \qbezier(65, 80)(115, 80)(124, 62)
\put(65, 77){\makebox(0,0)[ct]{$\Phi_d$}} \put(6, 62){\vector(-1, -2){1}}
\qbezier(5, 20)(15, 0)(65, 0) \qbezier(65, 0)(115, 0)(124, 18)
\put(65, 3){\makebox(0,0)[cb]{$\Phi_d$}} \put(6, 18){\vector(-1, 2){1}}
\end{picture}}\end{picture} \caption[$G$ and $H$]{\label{Fdisk}
Construction and straightening of $g_c$ by employing mappings in the
exterior of the unit disk. If the filled Julia sets are not connected, the
diagram is well-defined and commuting on smaller neighborhoods of $\infty$.}
\end{figure}

\subsection[Construction of the Quasi-Quadratic Mapping g]%
 {Construction of the Quasi-Quadratic Mapping $g_c$} \label{32}
For $c\in\E_\sM$, we construct a quasi-regular mapping $g_c$ coinciding with
$g_c^\pre$ on $\K_c\,$. By employing the Boettcher conjugation $\Phi_c\,$,
the work will be done in the exterior of the unit disk $\overline\disk$.
This is convenient when $c\in\E_\sM\,$, and essential to construct the
homeomorphism $h$ in the exterior. In $\hat\C\setminus\K_c$ we have
$g_c^\pre=\Phi_c^{-1}\circ G^\pre\circ\Phi_c\,$, where $G^\pre:\csd\to\csd$
is discontinuous on six straight rays, and given by compositions of
$F(z):=z^2$ in the regions between these rays --- it is independent of $c$
in particular.
\begin{enumerate}
 \item First construct smooth domains $U,\,U'$ with
$\overline\disk\subset U$ and $\overline U\subset U'\,$, and a smooth
mapping $G:U\setminus\overline\disk\to U'\setminus\overline\disk$. It shall
be proper of degree $2$ and coincide with $G^\pre$ except in sectors around
those six rays $\r(\Theta_i^\pm)$, where $G^\pre$ has a shift discontinuity.
The sectors are of the form $|\arg z-2\pi\Theta_i^\pm|<s\log|z|$, and there
$G$ is chosen conveniently as a $1$-homogeneous function of
 $\log z-\i2\pi\Theta_i^\pm$, thus ensuring that the dilatation bound does
not explode at the vertex of the sector. (This simplifies the construction
of \cite{bfe}, which employed a pullback of quadrilaterals.) The domains are
chosen in a finite recursion, employing that some iterate of $G^\pre$ is
strictly expanding \cite[Sect.~5.2]{wjt}. Since any orbit is visiting at
most two of the sectors, the dilatation of all iterates of $G$ is bounded
uniformly.
 \item Choose the radius $R>1$ and the conformal mapping
 $H:\hat\C\setminus\overline{U'}\to\hat\C\setminus\overline\disk_{R^2}$
with the normalization $H(z)=z+\O(1/z)$ at $\infty$ (which determines $R$
and $H$ uniquely). Extend $H$ to a quasi-conformal mapping
 $H:\hat\C\setminus U\to\hat\C\setminus\disk_R$
with $F\circ H=H\circ G$ on $\partial U$. Define the extended
$G:\csd\to\csd$ by $G:=H^{-1}\circ F\circ H$ on $\hat\C\setminus U$.
Now $G$ is proper of degree $2$, quasi-regular, and the dilatation of $G^n$
is bounded by some $K$ uniformly in $n$. Finally, extend $H$ to a mapping
$H:\csd\to\csd$ by recursive pullbacks, such that $F\circ H=H\circ G$
everywhere, then $H$ is $K$-quasi-conformal. Cf.~Fig.~\ref{Fdisk}.
 \item Now, set $g_c:=g_c^\pre$ on $\K_c$ and
 $g_c:=\Phi_c^{-1}\circ G\circ\Phi_c$ on $\hat\C\setminus\K_c\,$. Then $g_c$
is a quasi-quadratic mapping, i.e.~proper of degree $2$, with a uniform
bound on the dilatation of the iterates, with $\overline\partial g_c=0$
a.e.~on $\K_c\,$, and analytic in a neighborhood of $\infty$ with
$g_c(z)=z^2+\O(1)$. (It is continuous at $\gamma_c(\Theta_i^\pm)$ by
Lindel\"of's Theorem.)
\end{enumerate}

Now suppose that $c\in\P_\sM\setminus\E_\sM$ with $\Phi_\sM(c)\in U'$.
Then $\K_c$ is totally disconnected, and $\Phi_c$ is not defined in all of
$\hat\C\setminus\K_c\,$. It can be defined, however, in a domain mapped to
the six sectors and to $\hat\C\setminus\overline U$ by $\Phi_c$
 \cite{bfe, wjt}. Thus $g_c$ is defined in this case as well, by matching
$g_c^\pre$ with $\Phi_c^{-1}\circ G\circ\Phi_c\,$.

In the following section, we shall construct an invariant ellipse field for
the quasi-quadratic $g_c\,$, and employ it to straighten $g_c\,$, i.e.~to
conjugate it to a quadratic polynomial $f_d\,$. Then we set $h(c):=d$. If we
had skipped step~2, $g_c$ would not be a quasi-quadratic mapping $\C\to\C$,
but a quasi-regular quadratic-like mapping (cf.~\cite{dhp, wjt}) between
bounded domains $U_c\to U_c'\,$. This distinction is related to possible
alternative techniques:

\begin{rmk}[Alternative Techniques] \label{Ra}
1. The classical techniques would be as follows \cite{bd, bfl}: after the
quasi-regular quadratic-like mapping $g_c:U_c\to U_c'$ is constructed, it is
not extended to $\C$, but it is first conjugated to an analytic
quadratic-like mapping, employing an invariant ellipse field in $U_c'\,$.
Then the latter mapping is straightened to a polynomial by the Straightening
Theorem \cite{dhp}. With this approach, it will not be possible to extend
the homeomorphism $h$ to the exterior of $\M$.

2. Here we shall use the same techniques as in \cite{bfe}: having extended
$g_c$ to $\hat\C$, it will be easy to straighten. Instead of applying the
Straightening Theorem, its proof \cite{cg} was adapted into the construction
of $g_c\,$. This approach makes the extension of $h$ to the exterior of $\M$
possible. By applying this technique to the construction of $g_c$ and $h(c)$
for $c\in\E_\sM$ as well, the proofs of bijectivity, continuity, and of
landing properties (Sect.~\ref{41}) are simplified.

3. Alternatively, $g_c:U_c\to U_c'$ could be constructed as a quasi-regular
quadratic-like mapping on a bounded domain, and be straightened without
extending it to $\hat\C$ first, by incorporating the alternative proofs of
the Straightening Theorem according to~\cite{dhp}. This proof is more
involved, but it has the advantage that the mapping $H$ can be chosen more
freely on $\overline{U'}\setminus U\,$, e.g.~such that it is the identity
on $\r(\Theta_1^\pm)$ and $\r(\Theta_4^\pm)$ \cite{wjt}. Then $h$ would be
the identity on the corresponding parameter rays, which makes it easier
to paste different homeomorphisms together.
\end{rmk}

\subsection[h is a Homeomorphism]{$h$ is a Homeomorphism} \label{33}
For $c\in\E_\sM\,$, or $c\in\P_\sM\setminus\E_\sM$ with
 $\Phi_\sM(c)\in U'$, the quasi-quadratic mapping $g_c$ was constructed
in the previous section. Now construct the $g_c$-invariant ellipse field
$\mu$ by pullbacks with $g_c\,$, such that $\mu=0$ in a neighborhood of
$\infty$ and a.e.~on $\K_c\,$. It is bounded by $(K-1)/(K+1)$, since the
dilatation of all iterates $g_c^n$ is bounded by $K$. Denote by $\psi_c$ the
solution of the Beltrami equation $\overline\partial\psi=\mu\partial\psi$,
normalized by $\psi_c(z)=z+\O(1/z)$, which is mapping the infinitesimal
ellipses described by $\mu$ to circles. Now $\psi_c\circ
g_c\circ\psi_c^{-1}$ is analytic on $\hat\C$ and proper of degree $2$, thus
a quadratic polynomial of the form $f_d(z)=z^2+d$. In a neighborhood of
$\infty$, $H\circ\Phi_c\circ\psi_c^{-1}$ is conjugating $f_d$ to $F$,
cf.~~Fig.~\ref{Fdisk}. By the uniqueness of the Boettcher conjugation, this
mapping equals $\Phi_d\,$. Recursive pullbacks show equality in
$\hat\C\setminus\K_d\,$, if $\K_c$ and $\K_d$ are connected, i.e.~for
$c\in\E_\sM\,$. Otherwise, equality holds on an $f_d$-forward-invariant
domain of $\Phi_d\,$, which may be chosen to include the critical value $d$.

We set $h(c):=d=\psi_c(c)$. If $c\in\E_\sM\,$, a combinatorial argument
shows $d\in\E_\sM\,$. If $c\in\P_\sM\setminus\E_\sM$ with
 $\Phi_\sM(c)\in U'$, we have
\be\label{hPhiH} \Phi_\sM(d)=\Phi_d(d)=\Phi_d\circ\psi_c(c)
 =H\circ\Phi_c(c)=H\circ\Phi_\sM(c) \ . \ee
Denote by $\tilde\P_\sM$ the closed strip that is bounded by the four curves
$\Phi_\sM^{-1}\circ H(\r(\Theta_i^\pm))$, $i=1,\,4$, which are quasi-arcs.
Now $h$ is extended to $h:\P_\sM\to\tilde\P_\sM$ by setting
\be\label{hH} h:=\Phi_\sM^{-1}\circ H\circ\Phi_\sM \,:\,
 \P_\sM\setminus\E_\sM\to\tilde\P_\sM\setminus\E_\sM \ . \ee
By (\ref{hPhiH}), this agrees with the definition of $h(c)$ by straightening
$g_c\,$, if $\Phi_\sM(c)\in U'$. Now (\ref{hH}) shows that $h$ is
bijective and $K$-quasi-conformal in the exterior of $\E_\sM\,$. We will see
that $h$ is bijective and continuous on $\E_\sM\,$. Let us remark that for
$c\in\E_\sM\,$, the value of $d=h(c)$ does not depend on the choices made in
the construction of $G$ and $H$, since $\psi_c$ is a hybrid-equivalence
\cite{dhp}. The proof of bijectivity in \cite{bfl} relied on this
independence, but the following one is simplified by employing $H$:

For $d\in\E_\sM\,$, consider $\tilde g_d^\pre$ according to Def.~\ref{Dg},
and define the quasi-quadratic mapping $\tilde g_d$ with
 $\tilde g_d:=\tilde g_d^\pre$ on $\K_d$ and
 $\tilde g_d:=\Phi_d^{-1}\circ\tilde G\circ\Phi_d$ in
$\hat\C\setminus\K_d\,$, where $\tilde G:=H\circ F\circ H^{-1}$. To see that
this choice is possible, note that $H$ is mapping the region
 $V\subset U\setminus\overline\disk$ (corresponding to $V_c$) to a distorted
version of $\tilde V$. There we have
 $\tilde G^\pre=F^{2-k_v}\circ(\pm F^{\tilde k_v-1})$. Observing that
 $F=H\circ G\circ H^{-1}$ and $H$ commutes with $\pm\mathrm{id}$ on the set
in question, we have
 $\tilde G^\pre=H\circ G^{2-k_v}\circ(\pm G^{\tilde k_v-1})\circ H^{-1}$.
Following the orbit and applying the piecewise definition of $G^\pre$ yields
$\tilde G^\pre=H\circ F\circ H^{-1}$. Together with the same result in other
regions, this justifies the definition of $\tilde G$, i.e.~$\tilde g_d$ is
quasi-quadratic. Now $\tilde h(d)$ is defined by straightening
 $\tilde g_d\,$. --- Suppose that $c\in\E_\sM$ and $d=h(c)$, then
 $f_d=\psi_c\circ g_c\circ\psi_c^{-1}$ and by its definition in terms of
 $H=\Phi_d\circ\psi_c\circ\Phi_c^{-1}\,$, we have
 $\tilde g_d=\psi_c\circ f_c\circ\psi_c^{-1}$. Therefore $c=\tilde h(d)$
and $\tilde\psi_d=\psi_c^{-1}$. $\tilde h\circ h=\mathrm{id}$ and the
converse result imply that $h:\E_\sM\to\E_\sM$ is bijective with
$h^{-1}=\tilde h$.

By (\ref{hH}), $h$ is quasi-conformal in the exterior of $\E_\sM\,$. The
interior of $\E_\sM$ consists of a countable family of hyperbolic
components, plus possibly a countable family of non-hyperbolic components.
The former are parametrized by multiplier maps, the latter by transforming
invariant line fields. In both cases, $h$ is given by a composition of these
analytic parametrizations \cite{bfl, wjt}. It remains to show that $h$ is
continuous at $c_0\in\partial\E_\sM$: suppose $c_n\to c_0\,$, $d_n=h(c_n)$,
$d_0=h(c_0)$. By bijectivity we have
$d_0\in\partial\E_\sM=\E_\sM\cap\partial\M$. It does not matter if $c_n$
belongs to $\E_\sM$ or not. (Now we employ the definition of $h$ by
straightening $g_c\,$, which is equivalent to (\ref{hH}). One special case
requires extra treatment: when some $\gamma_{c}(\Theta_i^\pm)$ is iterated
to $\gamma_c(\Theta_1^\pm)$, and $c_0=\gamma_\sM(\Theta_i^\pm)$, then
$g_{c_0}$ is not defined.) It is sufficient to show
 $d_n\to d_*\Rightarrow d_*=d_0\,$. Since the $K$-quasi-conformal mappings
$\psi_n$ are normalized, there is a $K$-quasi-conformal $\Psi$ and a
subsequence $\psi_{c_n'}\to\Psi$, uniformly on $\hat\C$ \cite{lvqk}. We have
$\psi_{c_n'}\circ g_{c_n'}\circ\psi_{c_n'}^{-1}\to\Psi\circ
g_{c_0}\circ\Psi^{-1}$ and
 $\psi_{c_n}\circ g_{c_n}\circ\psi_{c_n}^{-1}=f_{d_n}\to f_{d_*}\,$, thus
$\Psi\circ\psi_{c_0}^{-1}$ is a quasi-conformal conjugation from $f_{d_0}$
to $f_{d_*}\,$. Although it need not be a hybrid-equivalence,
$d_0\in\partial\M$ implies $d_*=d_0$ \cite{dhp}. By the same arguments, or
by the Closed Graph Theorem, $h^{-1}$ is continuous as well. Thus
$h:\P_\sM\to\tilde\P_\sM$ is a homeomorphism mapping $\E_\sM\to\E_\sM\,$.

\subsection[Properties of h]{Further Properties of $h$} \label{34}
Since $h$ is analytic in the interior of $\E_\sM$ and quasi-conformal in the
exterior, it is natural to ask if it is quasi-conformal everywhere. Branner
and Lyubich are working on a proof employing quasi-regular quadratic-like
germs. Maybe an alternative proof can be given by constructing a homotopy
from $f_c$ to $\tilde g_d\,$, thus from $\mathrm{id}$ to $h$.

The dynamics of $h$ on $\E_\sM$ is simple: set
$c_0:=\gamma_\sM(\Theta_2^\pm)$ and $c_n:=h^n(c_0)$, $n\in\Z$. The connected
component of $\E_\sM$ between the two pinching points $c_n$ and $c_{n+1}$ is
a fundamental domain for $h^{\pm 1}$. These are accumulating at the
Misiurewicz points $a$ and $b$, and the method of \cite{tls} yields a linear
scaling behavior. Thus $h$ and $h^{-1}$ are Lipschitz continuous at $a$ and
$b$ (and H\"older continuous at all Misiurewicz points). For
$c\in\E_\sM\setminus\{a,\,b\}$ we have $h^n(c)\to b$ as $n\to\infty$ and
$h^n(c)\to a$ as $n\to-\infty$.

\section{Related Results and Possible Generalizations} \label{4}
Further results and examples from \cite{wjt} are sketched, and we present
some ideas on surgery for general one-parameter families.

\subsection{Combinatorial Surgery} \label{41}
The unit circle $\partial\disk$ is identified with $S^1:=\R/\Z$ by the
parametrization $\exp(\i2\pi\theta)$. For $h$ constructed from $g_c^\pre$
according to Thm.~\ref{Th}, recall the mappings $F,\,G,\,H:\csd\to\csd$ from
Sect.~\ref{32}. Denote their boundary values by $\bF,\,\bG,\,\bH:S^1\to
S^1$. Thus $\bF(\theta)=2\theta\mod1$ and $\bG$ is piecewise linear. Now
$\bH$ is the unique orientation-preserving circle homeomorphism conjugating
$\bG$ to $\bF$, $\bH\circ\bG\circ\bH^{-1}=\bF$. $\bH(\theta)$ is computed
numerically from the orbit of $\theta$ under $\bG$ as follows: for $n\in\N$,
the $n$-th binary digit of $\bH(\theta)$ is $0$ if
$0\le\bG^{n-1}(\theta)<1/2$, and $1$ if $1/2\le\bG^{n-1}(\theta)<1$. For
rational angles, the (pre-) periodic sequence of digits is obtained from a
finite algorithm.

In the exterior of $\E_\sM\,$, $h$ is represented by $H$ according to
(\ref{hH}). Applying this formula to parameter rays and employing
Lindel\"of's Theorem shows: $\r_\sM(\theta)$ is landing at
$c\in\partial\E_\sM$, iff $\r_\sM(\bH(\theta))$ is landing at $h(c)$. If $c$
is a Misiurewicz point or a root, then $\theta$ is rational, and
$\bH(\theta)$ is computed exactly. In this sense, $d=h(c)$ is determined
combinatorially. Alternatively, one can construct the critical orbit of
$g_c^\pre$ and the Hubbard tree of $f_d\,$. The simplest case is given when
the critical orbit meets $\E_c$ only once: then the orbit of $c$ under
$g_c^\pre$ is the same as the orbit of $\eta_c(c)$ under $f_c\,$.

Regularity properties of $\bH$ are discussed in \cite[Sect.~9.2]{wjt}. $\bH$
has Lipschitz or H\"older scaling properties at all rational angles. $\bH$
and $\bH^{-1}$ are H\"older continuous with the optimal exponents $\tilde
k_v/k_v$ and $k_w/\tilde k_w\,$. Since $H$ is $K$-quasi-conformal, Mori's
Theorem \cite{lvqk} says that $\bH^{\pm1}$ is $1/K$-H\"older continuous.
Thus we have the lower bound $K\ge\max(k_v/\tilde k_v,\,\tilde k_w/k_w)$,
independent of the choices made in the construction of
$h:\P_\sM\setminus\E_\sM\to\tilde\P_\sM\setminus\E_\sM\,$. By a piecewise
construction we obtain a homeomorphism $h:\M\to\M$, which extends to a
homeomorphism of $\C$, but such that no extension can be quasi-conformal.

\subsection{Homeomorphisms at Misiurewicz Points} \label{42}
A homeomorphism $h:\E_\sM\to\E_\sM$ according to Thm.~\ref{Th} is expanding
at the Misiurewicz point $a$. Asymptotically, $\M$ shows a linear scaling
behavior at $a$. (In Fig.~\ref{Fedge}, you can observe the asymptotic
self-similarity of $\M$ at $a$, and similarity between $\M$ at $c\approx
c_n$ and $\K_{c_n}$ at $z\approx0$.) Now $h$ is asymptotically linear in a
``macroscopic'' sense, e.g.~there is an asymptotically linear sequence of
fundamental domains, but this is not true pointwise. These results are
obtained by combining the techniques from \cite{tls} with the combinatorial
description of $h$ according to Sect.~\ref{41}: consider a suitable sequence
$c_n\to a$. If the critical orbit of $f_{c_n}$ travels through $\E_{c_n}$
once, then $h$ is asymptotically linear on the sequence, but it is not if
the orbit meets $\E_{c_n}$ twice.

Conversely, given a branch at some Misiurewicz point $a$, is there an
appropriate homeomorphism $h$? We only need to find a combinatorial
construction of $g_c^\pre\,$. This is done in \cite{wjt} e.g.~for all
$\beta$-type Misiurewicz points. (Here $\P_\sM$ and $\overline{V_c\cup W_c}$
are sectors, not strips.) The author's research was motivated by discussions
with D.~Schleicher, who had worked on the construction of dynamics in the
parameter plane before.

\subsection{Edges, Frames, and Piecewise Constructions} \label{43}
For parameters $c$ in the $p/q$-limb of $\M$, the filled Julia set $\K_c$
has $q$ branches at the fixed point $\alpha_c$ of $f_c\,$. A connected
subset $\E_c\subset\K_c$ is a \emph{dynamic edge} of order $n$, if
$f_c^{n-1}$ is injective on $\E_c$ and $f_c^{n-1}(\E_c)$ is the part of
$\K_c$ between $\alpha_c$ and $-\alpha_c\,$. (More precisely, $f_c^{n-1}$
shall be injective in a neighborhood of the edge without its vertices.) The
edge is characterized by the external angles at the vertices. As $c$ varies,
it may still be defined, or it may cease to exist after a bifurcation of
preimages of $\alpha_c\,$. Now $\E_\sM\subset\M_{p/q}$ is a \emph{parameter
edge}, if for all $c\in\E_\sM$ the dynamic edge $\E_c$ (with given angles)
exists and satisfies $c\in\E_c\,$, and if $\E_\sM$ has the same external
angles as $\E_c\,$. In Figs.~\ref{Fstrips}--\ref{Fedge}, $\E_\sM$ is the
parameter edge of order $4$ in $\M_{1/3}\,$.

$\M_{p/q}$ contains a little Mandelbrot set $\M'=c_0\ast\M$ of period $q$
(cf.~Sect.~\ref{44}). If a parameter edge $\E_\sM$ is behind $c_0\ast(-1)$,
there is a homeomorphism $h:\E_\sM\to\E_\sM$ analogous to that of
Figs.~\ref{Fstrips}--\ref{Fedge}. Behind the $\alpha$-type Misiurewicz point
$c_0\ast(-2)$, edges can be decomposed into subedges and \emph{frames}
\cite[Sect.~7]{wjt}. These frames are constructed recursively, like the
intervals in the complement of the middle-third Cantor set. A family of
homeomorphisms on subedges shows that all frames on the same edge are
mutually homeomorphic, and they form a finer decomposition than the
fundamental domains of a single homeomorphism. By permuting the frames (in a
monotonous way), new homeomorphisms $h$ are defined piecewise. These may
have properties that are not possible when $h$ is constructed from a single
surgery. E.g., in contradiction to Sect.~\ref{34}, $h$ can be constructed
such that it is not Lipschitz continuous or not even H\"older continuous at
the vertex $a$ of $\E_\sM\,$. Or it can map a Misiurewicz point with two
external angles to a parameter with irrational angles, which is not a
Misiurewicz point.

The notions of edges and frames can be generalized: for parameters $c$
behind the root of a hyperbolic component, the filled Julia set $\K_c$
contains two corresponding pre-characteristic points, which take the roles
of $\pm\alpha_c\,$.

\subsection{Tuning and Composition of Homeomorphisms} \label{44}
\begin{figure}[h!t!b!]
\unitlength 1mm \linethickness{0.4pt}
\begin{picture}(150, 45)
\put(5, 0){\includegraphics[width=60mm]{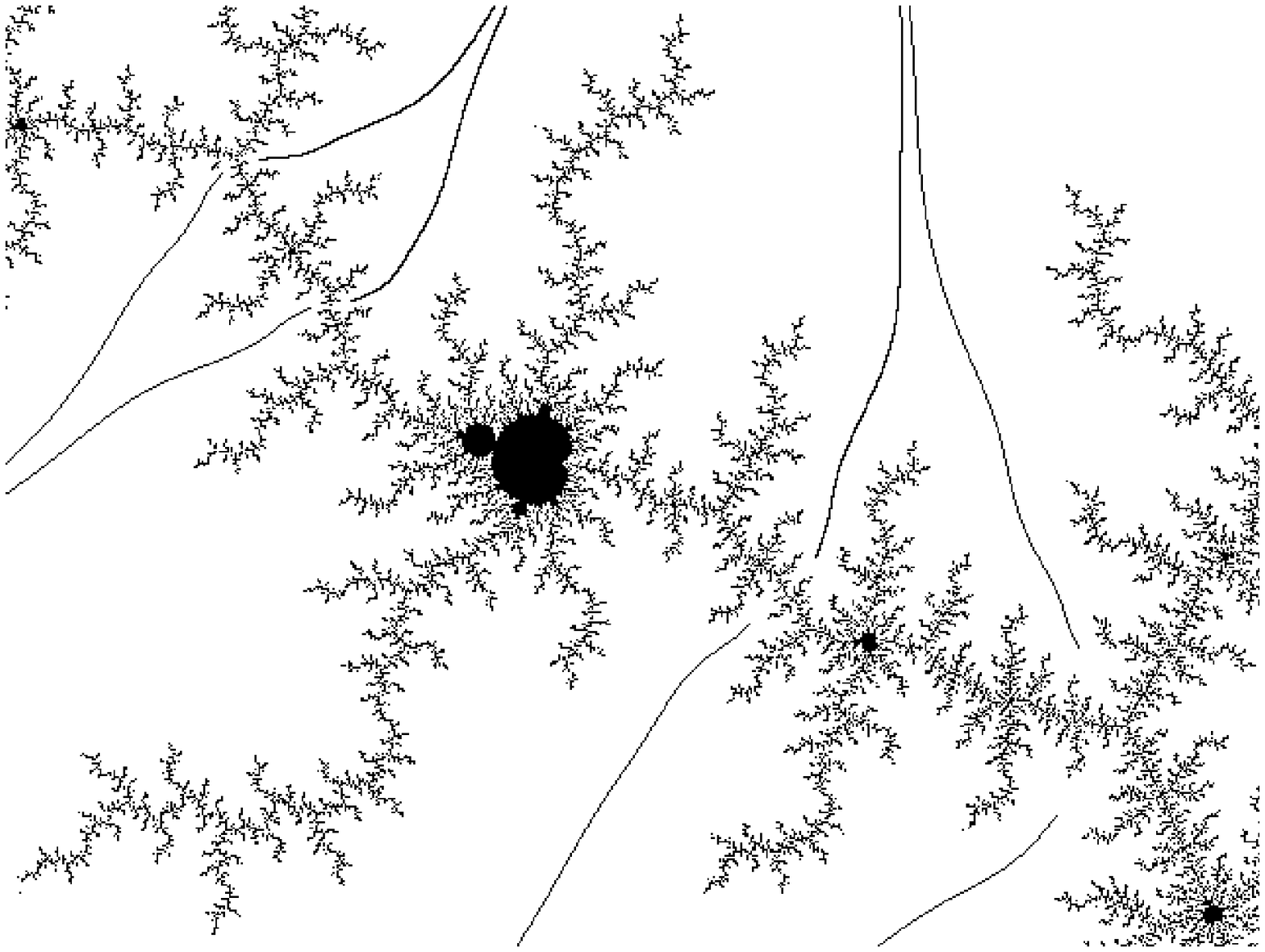}}
\put(85, 0){\includegraphics[width=60mm]{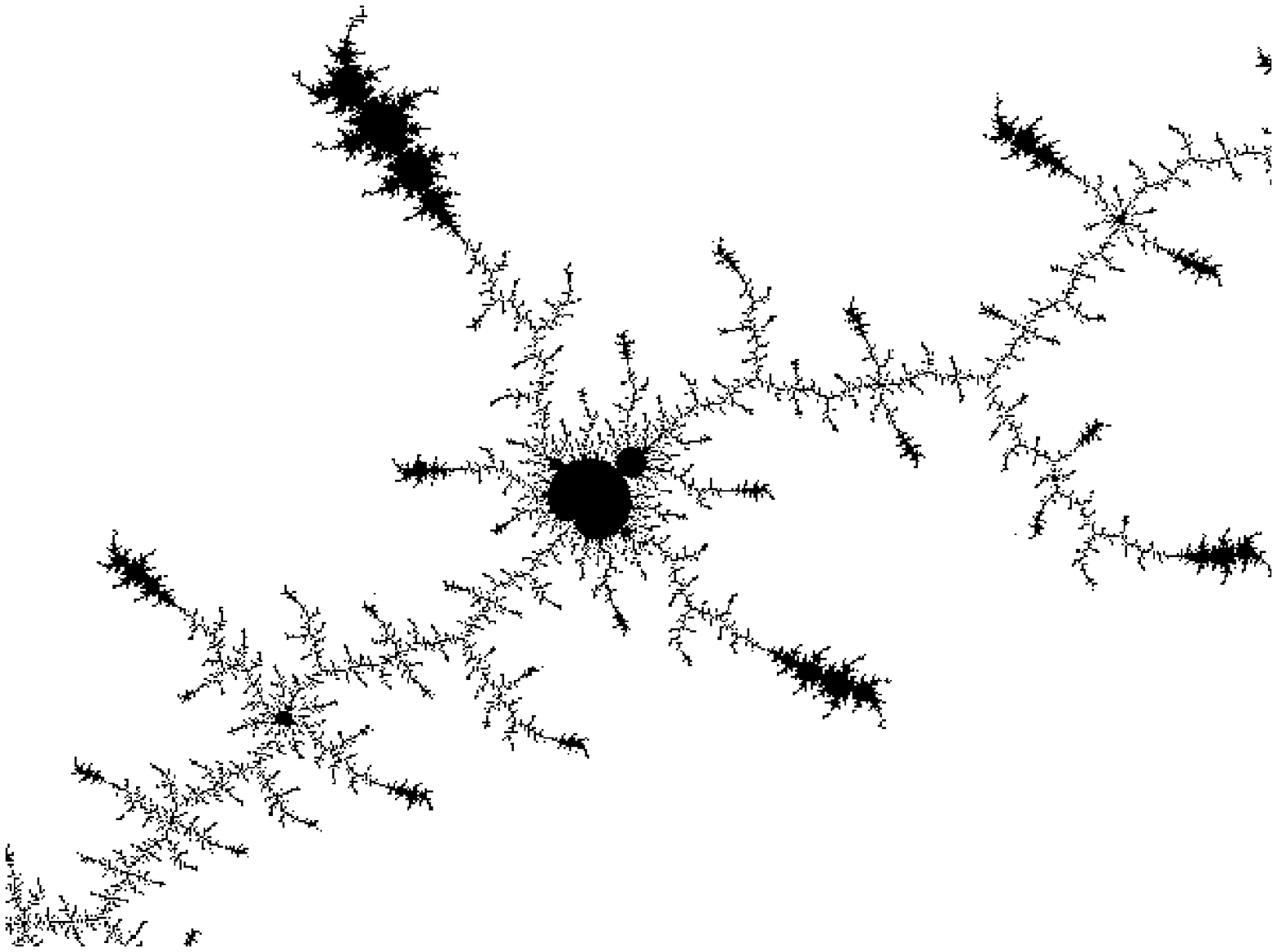}}
\thicklines
\put(66, 11){\makebox(0,0)[cc]{$a$}}
\put(10, 34){\makebox(0,0)[cc]{$b$}}
\put(48, 27){\vector(-1, 1){11}}
\put(44.5, 37){\makebox(0,0)[cc]{$h'$}}
\put(86, 8){\makebox(0,0)[cc]{$a$}}
\put(130 ,32){\makebox(0,0)[cc]{$b$}}
\put(93, 21){\vector(2, 3){9}}
\put(94, 29.5){\makebox(0,0)[cc]{$h'$}}
\end{picture} \caption[]{\label{Ftusi}
Two homeomorphisms $h':\E_\sM'\to\E_\sM'$ obtained from a similar
construction as $h:\E_\sM\to\E_\sM$ in Figs.~\ref{Fstrips}--\ref{Fedge}.
Left: tuning with the center $c_0=-1$ yields an edge in the limb
$\M_{1/2}\subset\M$. The eight angles $\Theta_i^\pm$ are obtained by tuning
those of Fig.~\ref{Fstrips}, i.e.~replacing the digits $0$ by $01$ and $1$
by $10$. $h'$ is defined not only on $c_0\ast\E_\sM\,$, but on a strip
including all decorations. Right: part of the parameter space of cubic
polynomials with a persistently indifferent fixed point. The connectedness
locus contains copies of a quadratic Siegel Julia set \cite{bhjps}. Again,
$h'$ is defined in a strip containing a countable family of decorations
attached to the copy of $\E_\sM\,$.}
\end{figure}

For a center $c_0$ of period $p$, there is a ``little Mandelbrot set''
$\M'\subset\M$ and a tuning map $\M\to\M'$, $x\mapsto y=c_0\ast x$ with
$0\mapsto c_0\,$. Now $\K_y$ contains a ``little Julia set'' $\K_{y,\,p}$
around $0$, where $f_y^p$ is conjugate to $f_x$ on $\K_x\,$
 \cite{dhp, TLHa}. A homeomorphism $h:\E_\sM\to\E_\sM$ according to
Thm.~\ref{Th} is compatible with tuning in two different ways:
\begin{itemize}
 \item If $c_0\in\E_\sM\,$, then $h$ is mapping $\M'$ to the little
Mandelbrot set at $h(c_0)$:\\
 $h(c_0\ast x)=(h(c_0))\ast x$. Cf.~\cite{bfl}.
 \item For any center $c_0\in\M$, set $\E_\sM':=c_0\ast\E_\sM\subset\M$. A
new homeomorphism $h':\E_\sM'\to\E_\sM'$ is obtained by composition,
i.e.~$h'(c_0\ast x):=c_0\ast(h(x))$. Now $\E_\sM'$ is obtained by
disconnecting $\M$ at a countable family of pinching points, but $h'$ has a
natural extension to all of these ``decorations'' (except for two): the
mapping $\eta_c$ that produced the homeomorphism $h$ is transferred by
cutting the little Julia set into strips. The required pinching points do
not bifurcate when the parameter $y$ is in a decoration of $\E_\sM'\,$, thus
the new piecewise construction $\eta_y'$ works in a whole strip. An example
is shown in Fig.~\ref{Ftusi} (left).
\end{itemize}

The same principle applies, e.g., to crossed renormalization \cite{rscr}, or
to the Branner-Douady homeomorphism
$\Phi_A:\M_{1/2}\to\mathcal{T}\subset\M_{1/3}\,$: suppose that
$\E_\sM\subset\M_{1/2}$ and $h:\E_\sM\to\E_\sM$ is constructed according to
Thm.~\ref{Th}, i.e.~from a combinatorial $g_c^\pre$ according to
Def.~\ref{Dg}. Then $\E_\sM':=\Phi_A(\E_\sM)$ is a subset of $\M_{1/3}$,
where a countable family of decorations was cut off. Again
 $h':=\Phi_A\circ h\circ\Phi_A^{-1}:\E_\sM'\to\E_\sM'$ extends to a whole
strip by transferring the combinatorial construction of $g_c^\pre\,$. If,
e.g., $h$ is a suitable homeomorphism on the edge from $\gamma_\sM(5/12)$
to $\gamma_\sM(11/24)$, then $h'$ is the homeomorphism of Figs.~\ref{Fstrips},
\ref{Fedge}.

\subsection{Other Parameter Spaces} \label{45}
In Thm.~\ref{Th} we obtained homeomorphisms $h:\E_\sM\to\E_\sM$ of suitable
subsets $\E_\sM\subset\M$, but the method is not limited to quadratic
polynomials. To apply it to other one-dimensional families of polynomials or
rational mappings, these mappings must be characterized dynamically. The
polynomials of degree $d$ form a $(d-1)$-dimensional family (modulo affine
conjugation). Suppose that a one-dimensional subfamily $f_c$ is defined by
one or more of the following \emph{critical relations}:
\begin{itemize}
 \item A critical point of $f_c$ is degenerate, or one critical point is
iterated to another one, or critical orbits are related by $f_c$ being
even or odd.
 \item A critical point is preperiodic or periodic (superattracting).
 \item There is a persistent cycle with multiplier $\rho$,
$0<|\rho|\le1$. This cycle is always ``catching'' one of the critical
points, but the choice may change.
\end{itemize}
An appropriate combination of such relations defines a one-parameter family
$f_c\,$, where the coefficients and the critical points of $f_c$ are
algebraic in $c$. Locally in the parameter space, there is one
\emph{active} or \emph{free} critical point $\omega_c\,$, whose orbit
determines the qualitative dynamics. The other critical points are either
linked to $\omega_c\,$, or their behavior is independent of $c$. The
connectedness locus $\M_f$ contains the parameters $c$, such that the filled
Julia set $\K_c$ of $f_c$ is connected, equivalently
$f_c^n(\omega_c)\not\to\infty$, or $\omega_c\in\K_c\,$. In general
$\omega_c$ is not defined globally by an analytic function of $c$, since it
may be a multi-valued algebraic function of $c$, or since a persistent cycle
may catch different critical points. But looking at specific families, it
will be possible to define $\Phi_\sM$ and parameter rays for suitable
subsets of the parameter space, and to understand their landing properties.

Then an analogue of Thm.~\ref{Th} can be proved: for a piecewise defined
$g_c^\pre\,$, a quasi-polynomial mapping $g_c$ is constructed analogously to
Sect.~\ref{32}, and straightened to a polynomial. By the critical relations
and by normalizing conditions, it will be of the form $f_d\,$, and we set
$h(c):=d$. (At worst, the normalizing conditions will allow finitely many
choices for $d$.) Note that this procedure will not work, if our family
$f_c$ is an arbitrary submanifold of the $(d-1)$-dimensional family of all
polynomials, and not defined by critical relations.

When such a theorem is proved, the remaining creative step is the
combinatorial definition of $\E_\sM$ and $g_c^\pre$. Some examples can be
obtained in the following way: when a non-degenerate critical point is
active, the connectedness locus $\M_f$ will contain copies $\M'$ of $\M$
\cite{mmu}. Starting from a homeomorphism $h:\E_\sM\to\E_\sM\,$, $\M'$
contains a decorated copy of $\E_\sM\,$, and the corresponding homeomorphism
$h'$ extends to all decorations by an appropriate definition of $g_c^\pre$
--- the angles $\Theta_i^\pm$ are seen at the copy of a quadratic Julia set
within $\K_c\,$, where some iterate of $f_c$ is conjugate to a quadratic
polynomial. It remains to check that no other critical orbit is passing
through $\overline{V_c\cup W_c}\,$, then $g_c^\pre$ is well-defined. An
example is given in Fig.~\ref{Ftusi} (right).

The rational mappings of degree $d$ form a $(2d-2)$-dimensional family
(modulo M\"obius conjugation). Suppose that a one-dimensional subfamily
$f_c$ is defined by critical relations. When there are one or more
persistently (super-) attracting cycles, then $\K_c$ shall be the complement
of the basin of attraction, and $\M_f$ shall contain those parameters, such
that the local free critical point is not attracted. $\partial\M_f$ will be
the bifurcation locus \cite{mmu}. If the persistent cycles are
superattracting, we can define dynamic rays and parameter rays by the
Boettcher conjugation. When the topology and the landing properties are
understood sufficiently well, homeomorphisms can be constructed by
quasi-conformal surgery.

\begin{figure}[h!t!b!]
\unitlength 1mm \linethickness{0.4pt}
\begin{picture}(150, 63)
\put(0, 0){\includegraphics[width=84mm]{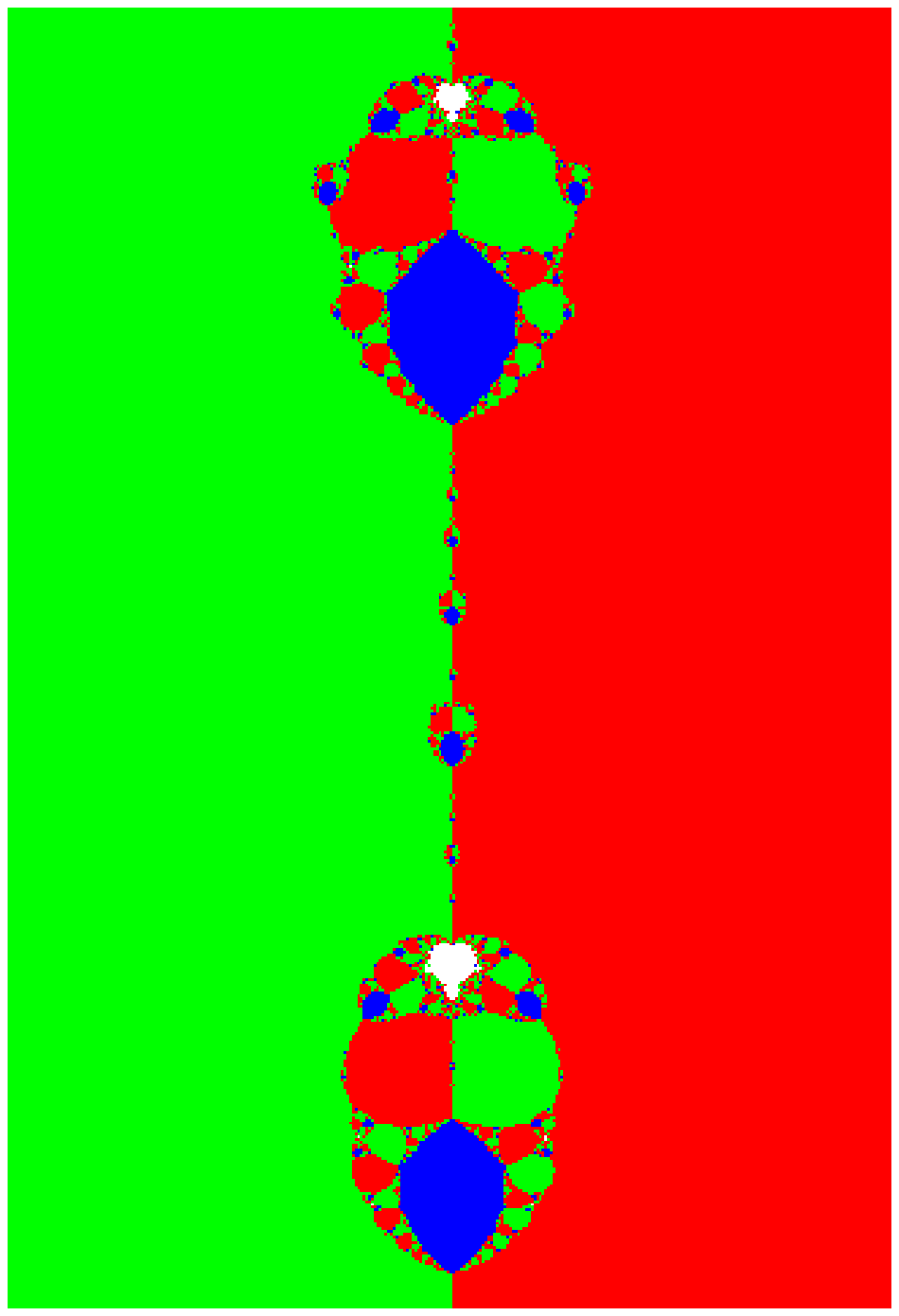}}
\put(66, 0){\includegraphics[width=84mm]{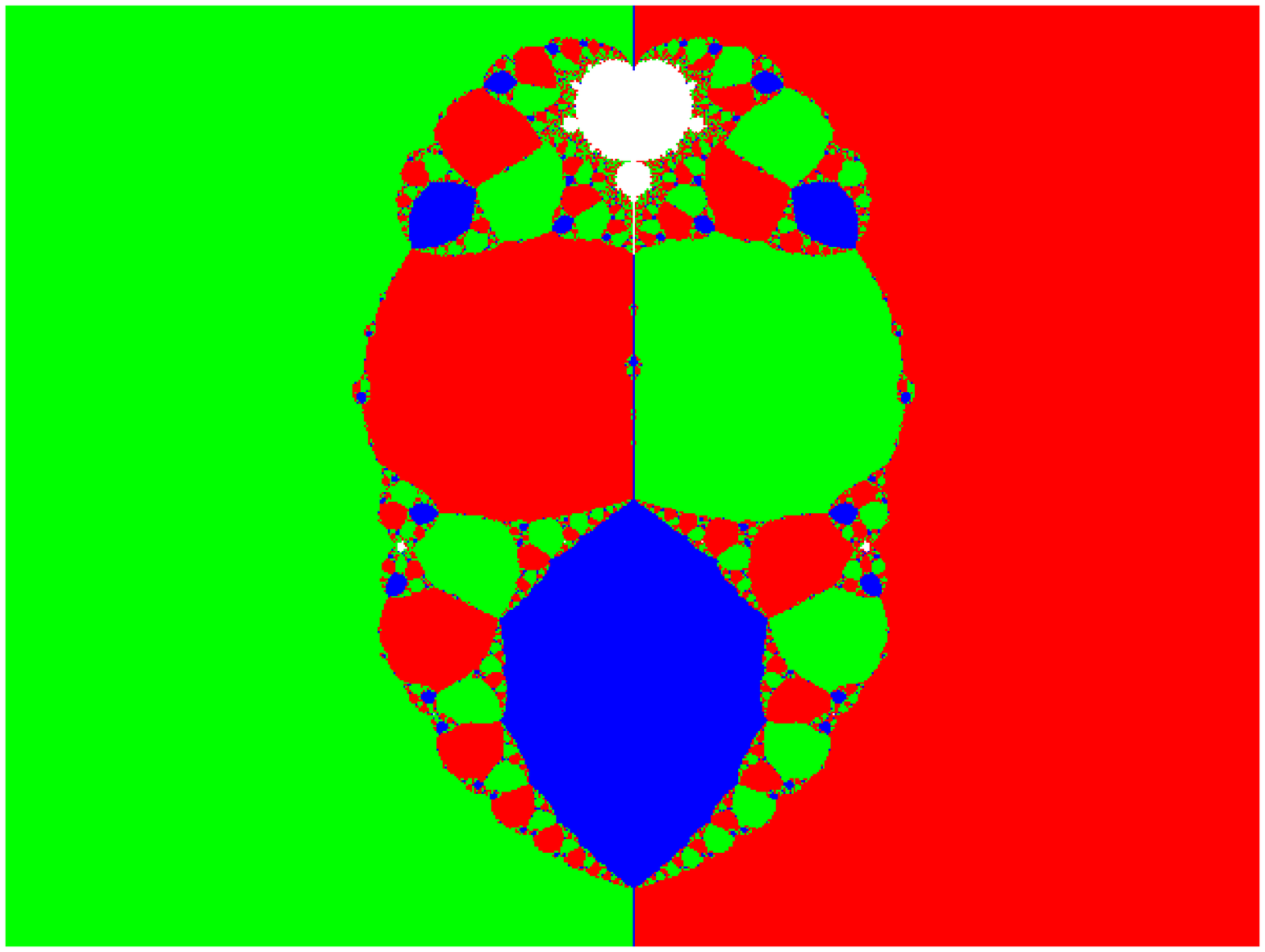}}
\thicklines
\put(18, 37){\vector(0, -1){10}}
\put(14, 33){\makebox(0,0)[cc]{$h_0$}}
\put(88, 45){\vector(-1, -4){2}}
\put(83, 43){\makebox(0,0)[cc]{$h_1$}}
\put(111, 32){\vector(0, 1){8}}
\put(116, 36){\makebox(0,0)[cc]{$h_2$}}
\end{picture} \caption[Newton method]{\label{Fn}
Homeomorphisms in the parameter space of Newton methods for cubic
polynomials. Left: an ``edge'' between the ``almonds'' of orders $3$ and
$2$. Right: homeomorphisms on edges within the almond of order $2$. (The
different colors, or shades of gray, indicate that the free critical point
is attracted to one of the roots of the corresponding polynomial.)}
\end{figure}

\begin{figure}[h!t!b!]
\unitlength 1mm \linethickness{0.4pt}
\begin{picture}(150, 45)
\put(0, 0){\includegraphics[width=60mm]{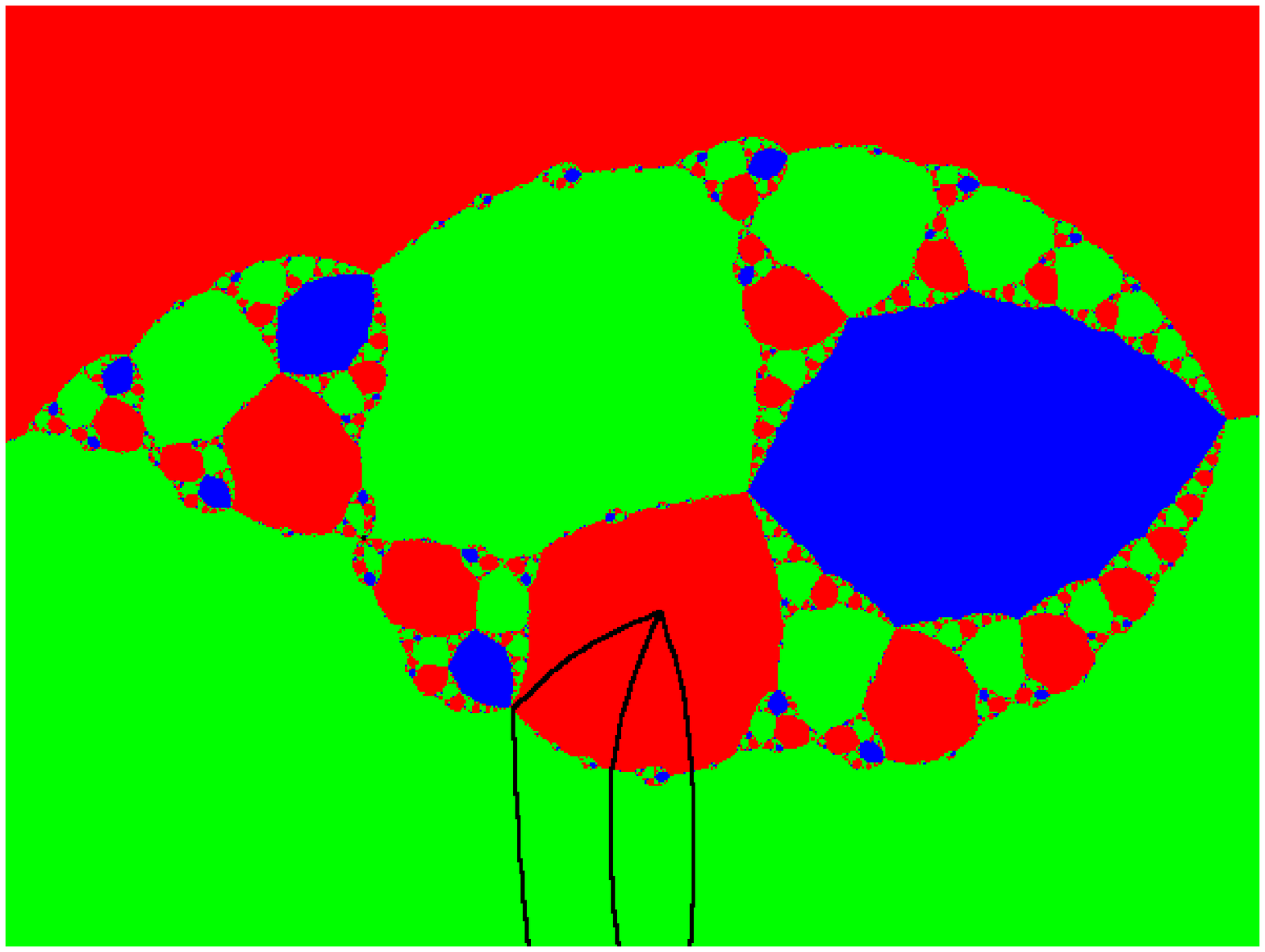}}
\put(90, 0){\includegraphics[width=60mm]{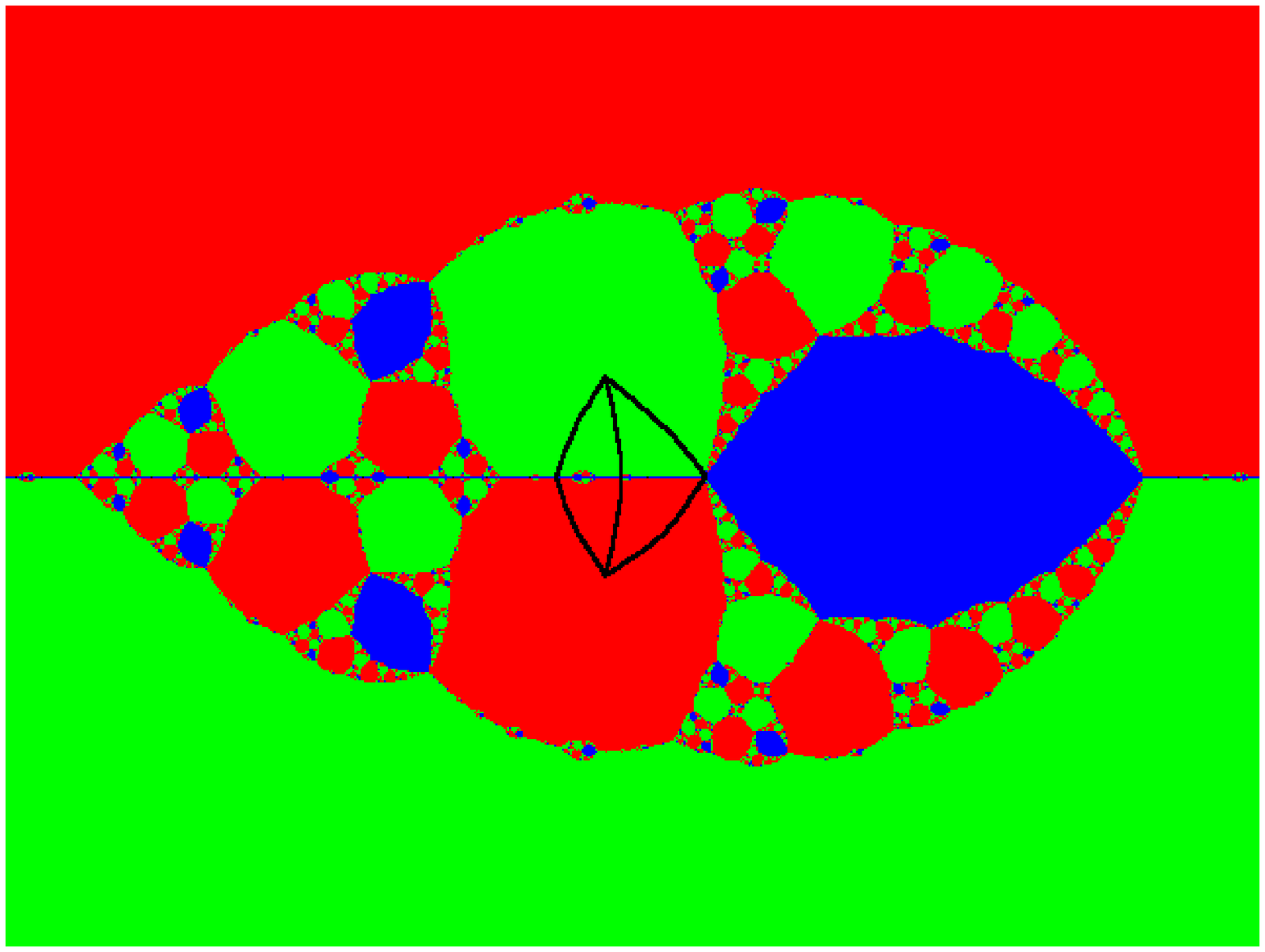}}
\end{picture} \caption[Newton method,dynamics]{\label{Fnj}
Cutting the Julia set with dynamic rays belonging to two different basins,
to define the strips $V_c\,,\,W_c$ and the mapping $g_c^\pre\,$. This yields
the homeomorphisms $h_1$ (left) and $h_2$ (right) in the almond of order $2$
(cf.~Fig.~\ref{Fn}).}
\end{figure}

An example is provided by cubic Newton methods, cf.~\cite{dhp, tlbcn, prt}:
$f_c$ has three superattracting fixed points and one free critical point.
Parts of the parameter space are shown in Fig.~\ref{Fn}. By dynamic rays in
the adjacent immediate basins of two fixed points, the Julia set is cut into
``strips'' to define $g_c^\pre$. In both basins, the
techniques of Sect.~\ref{32} are applied to construct the quasi-Newton
mapping $g_c\,$. It is straightened to $f_d\,$, and a homeomorphism is
obtained by $h_0(c):=d$. It is permuting little ``almonds,'' respecting
their decomposition into four colors. Similar constructions are possible
when one or both of the adjacent components of basins at $\E_c$ are not
immediate, i.e.~when the hyperbolic components at $\E_\sM$ are of greater
depth \cite{prt}. Cf.~$h_1\,,\,h_2$ in Figs.~\ref{Fn}--\ref{Fnj}.

Cubic Newton methods are understood as matings of cubic polynomials
\cite{tlbcn}, and there are analogous homeomorphisms in the parameter space
of cubic polynomials with one superattracting fixed point. Again, the rays
used in the piecewise definition of $g_c^\pre$ belong to the
basins of two attracting fixed points, but one is finite and one at $\infty$
in the polynomial case. H.~Hubbard has suggested to look at quadratic
rational mappings with a superattracting cycle, which contain matings of
quadratic polynomials. When we try to transfer a known homeomorphism of the
Mandelbrot set to this family, in general we will have to use articulated
rays to cut the Julia set. Although it may be possible to define $g_c^\pre$,
it will not be possible to construct the quasi-regular mapping $g_c\,$,
because the shift discontinuity happens not only within the basin of
attraction, but at pinching points of the Julia set as well. For the same
reason, it will not be possible to transfer a homeomorphism of $\M$ to a
neighborhood of a copy of $\M$ in the cubic Newton family.

\section[Homeomorphism Groups of M]{Homeomorphism Groups of $\M$} \label{5}
Denote the group of orientation-preserving homeomorphisms $h:\M\to\M$ by
$\G_\sM\,$. If two homeomorphisms coincide on $\partial\M$, they encode the
same information on the topological structure of $\M$. To avoid
trivialities, we suggest some definitions of groups of non-trivial
homeomorphisms as well.

\begin{dfn}[Groups of Homeomorphisms] \label{DG}
$1.$ $\G_\sM$ is the group of homeomorphisms $h:\M\to\M$ that are
orientation-preserving at branch points, and orientation-preserving in the
interior of $\M$.

$2.$ $\G_a$ is the group of homeomorphisms $h:\M\to\M$ that are
orientation-preserving at branch points, and analytic in the interior of
$\M$.

$3.$ $\G_b$ is the group of homeomorphisms $h:\partial\M\to\partial\M$ that
are orientation-preserving at branch points, and orientation-preserving on
the boundaries of hyperbolic components.

$4.$ $\G_q$ is the factor group $\G_\sM/\G_1\,$, where $\G_1$ is the normal
subgroup consisting of trivial homeomorphisms:
$\G_1:=\{h\in\G_\sM\,|\,h=\mathrm{id}\,\,\mbox{on}\,\,\partial\M\}$.
\end{dfn}

$\G_q$ is the most natural definition of non-trivial homeomorphisms.
$\G_a\,,\,\G_b\,,\,\G_q$ may well turn out to be mutually isomorphic. On
$\G_\sM\,,\,\G_a\,$, and $\G_b\,$, define a metric by
\ban\label{met} d(h_1,\,h_2) &:=&
 \|h_1-h_2\|_\infty + \|h_1^{-1}-h_2^{-1}\|_\infty \\[1mm]\nonumber
 &:=& \max|h_1(c)-h_2(c)| + \max|h_1^{-1}(c)-h_2^{-1}(c)| \ , \ean
where the maxima are taken over $c\in\M$ or $c\in\partial\M$, respectively.
$\G_q$ consists of equivalence classes of homeomorphisms coinciding on the
boundary, $[h]=h\G_1=\G_1h$. Since $\G_1$ is closed, a metric is given by
\ban d([h_1],\,[h_2]) &:=& \inf\,\Big\{\,  \nonumber
 \|h_1'-h_2'\|_\infty \,\Big|\, h_1'\in[h_1],\,h_2'\in[h_2] \,\Big\} \\[1mm]
 &+&  \inf\,\Big\{\, \|{h_1'}^{-1}-{h_2'}^{-1}\|_\infty \,\Big|\,
 h_1'\in[h_1],\,h_2'\in[h_2] \,\Big\}  \label{metq1} \\[1mm]
 &=& \inf\,\Big\{\, \|h_1\circ u-h_2\|_\infty +  \label{metq2}
 \|h_1^{-1}\circ v-h_2^{-1}\|_\infty \,\Big|\,u,\,v\in\G_1 \,\Big\} \ . \ean
It may be more natural to take the infimum of a sum instead of the sum of
infima in (\ref{metq1}), i.e.~$\inf d(h_1',\,h_2')$, but I do not know how
to prove the triangle inequality in that case. (\ref{metq2}) is obtained
from (\ref{metq1}) by employing the facts that $\G_1$ is normal, and that
right translations are isometries of the norm:
 $\|h_1-h_2\|_\infty=\|h_1\circ h-h_2\circ h\|_\infty\,$, since $h\in\G_\sM$
is bijective.

\begin{prop}[Topology of Homeomorphisms Groups] \label{PtG}
$1$. $\G_\sM\,,\,\G_a\,,\,\G_b\,,\,\G_q$ are complete metric spaces and
topological groups, i.e.~composition and inversion are continuous.

$2.$ For $\G=\G_\sM\,,\,\G_a\,,\,\G_b\,,\,\G_q$ we have: if
$\Omega_1\,,\,\Omega_2$ are hyperbolic components, then
 $\mathcal{N}:=\{h\in\G\,|\,h(\partial\Omega_1)=\partial\Omega_2\}$ is open.
\end{prop}

\textbf{Proof}: 1. For $\G_\sM\,,\,\G_a\,,\,\G_b\,$, the proof is
straightforward. But suppose we had used the alternative metric
 $\tilde d(h_1,\,h_2):=\|h_1-h_2\|_\infty\,$, and $(h_n)\subset\G_\sM$ is
a Cauchy sequence in that metric. Then it is converging uniformly to a
continuous, surjective $h:\E_\sM\to\E_\sM\,$. If $h$ is injective, then
$h_n^{-1}\to h^{-1}$ uniformly. But $h$ need not be injective, a
counterexample is constructed in item~2 of \cite[Prop.~7.7]{wjt}. Thus, if
we had used $\tilde d$ instead of $d$, the topology of
 $\G_\sM\,,\,\G_a\,,\,\G_b$ would be the same, but they would be incomplete
metric spaces.

Now suppose $([h_n])$ is a Cauchy sequence in $\G_q\,$. It is sufficient to
show that a subsequence converges, and without restriction we have
 $d([h_{n+1}],\,[h_n])\le3^{-n}$. Choose $u_n\,,\,v_n\in\G_1$ with
\[ \|h_{n+1}\circ u_n-h_n\|_\infty \le2^{-n} \quad\mbox{and}\quad
 \|h_{n+1}^{-1}\circ v_n-h_n^{-1}\|_\infty \le2^{-n} \ . \]
Define the sequences
\[ \hat h_n:= h_n\circ u_{n-1}\circ u_{n-2}\circ\dots\circ u_1
 \quad\mbox{and}\quad
 \tilde h_n:= h_n^{-1}\circ v_{n-1}\circ v_{n-2}\circ\dots\circ v_1 \ . \]
Since the maximum norm is invariant under right translations on $\M$, they
satisfy
\[ \|\hat h_{n+1}-\hat h_n\|_\infty \le2^{-n} \quad\mbox{and}\quad
 \|\tilde h_{n+1}-\tilde h_n\|_\infty \le2^{-n} \ , \]
and there are continuous functions $\hat h,\,\tilde h$ with
 $\hat h_n\to\hat h$ and $\tilde h_n\to\tilde h$ uniformly. Now
 $\hat h\circ\tilde h$ and $\tilde h\circ\hat h$ are uniform limits of a
sequence in $\G_1\,$, thus surjective, and $\hat h$ is a homeomorphism. We
have $\hat h_n\to\hat h$ in $\G_\sM$ and $[h_n]=[\hat h_n]\to[\hat h]$ in
$\G_q\,$, therefore $\G_q$ is complete.

2. Hyperbolic components can be distinguished topologically from
non-hyperbolic components, since only the boundary of a hyperbolic component
contains a countable dense set of pinching points (by the Branch Theorem
\cite{sf2}). Thus every homeomorphism of $\M$ or $\partial\M$ is permuting
the set of hyperbolic components or of their boundaries, respectively. Fix
$a,\,b\in\partial\Omega_2$, and choose $\eps>0$ such that no hyperbolic
component $\neq\Omega_2$ is meeting both of the disks of radius $\eps$
around $a$ and $b$. This is possible, since there are several external rays
landing at $\partial\Omega_2\,$. If $h_0\in\mathcal{N}$ and $h\in\G$ with
$d(h,\,h_0)<\eps$, then $|h(h_0^{-1}(a))-a|<\eps$ and
$|h(h_0^{-1}(b))-b|<\eps$, thus $h(\partial\Omega_1)=\partial\Omega_2\,$.
(Analogously for the classes in $\G_q\,$.) \mybox

\begin{thm}[Groups of Non-Trivial Homeomorphisms] \label{TG}
The groups of non-trivial homeomorphisms of $\M$ or $\partial\M$ ---
$\G_a\,$, $\G_b$ and $\G_q$ --- share the following properties:

$1.$ They have the cardinality of the continuum $\R$, and they are totally
disconnected.

$2.$ They are perfect, and not compact \(not even locally compact\).

$3.$ A family of homeomorphisms $\mathcal{F}\subset\G_a\,,\,\G_b\,,\,\G_q$
is called \emph{normal}, if its closure is sequentially compact. A necessary
condition is that for every hyperbolic component $\Omega\subset\M$, there
are only finitely many components of the form $h^{\pm1}(\Omega)$,
$h\in\mathcal F$. If $\M$ is locally connected, this condition will be
sufficient for $\mathcal F$ being normal.
\end{thm}

By composition, the homeomorphisms constructed by surgery according to
Thm.~\ref{Th} generate a countable subgroup of $\G_a\,$, $\G_b$ or $\G_q\,$.
Will it be dense? --- For $\G_\sM\,$, items~1 and~3 are wrong, and item~2 is
true but trivial. Hence the motivation to consider the groups of non-trivial
homeomorphisms. The same results hold for the analogous groups, where the
condition of preserving the orientation is dropped.

\textbf{Proof:} We prove the statements for $\G_a\,$, the case of $\G_b$ or
$\G_q$ is similar. There is a sequence of disjoint subsets $\E_n\subset\M$
with $\diam(\E_n)\to0$, and a sequence of analytic homeomorphisms
$h_n:\M\to\M$, such that $h_n=\mathrm{id}$ on $\M\setminus\E_n\,$,
$h_n\neq\mathrm{id}$. To construct these, fix a homeomorphism
$h_*:\E_\sM\to\E_\sM$ according to Thm.~\ref{Th}, e.g.~that of
Figs.~\ref{Fstrips} and \ref{Fedge}. Choose $\E_0\subset\E_\sM$ and a
homeomorphism $h_0:\E_0\to\E_0\,$, $h_0\neq\mathrm{id}$, such that $\E_0$ is
contained in a fundamental domain of $h_*\,$. This is possible e.g.~by the
tuning construction from Sect.~\ref{44}. Then set
 $h_n:=h_*^n\circ h_0\circ h_*^{-n}$ on $\E_n:=h_*^n(\E_0)$, and extend it
by the identity to a homeomorphism of $\M$. We have $\diam(\E_n)\to0$ by the
scaling properties of $\M$ at Misiurewicz points \cite{tls}. --- An
alternative approach is as follows: construct homeomorphisms
$h_n:\E_n\to\E_n\,$, such that $\E_n$ is contained in the limb $\M_{1/n}\,$,
then $\diam(\E_n)\to0$ by the Yoccoz inequality. These homeomorphisms can be
constructed by tuning, or at $\beta$-type Misiurewicz points according to
Sect.~\ref{42}, or on edges (Sect.~\ref{43}). All of the homeomorphisms
constructed below extend to homeomorphisms of $\C$, cf.~item~3 of
Remark~\ref{Ra}. (If $\M$ is locally connected, all homeomorphisms in
$\G_\sM\,$, $\G_a\,$ or $\G_b$ extend to homeomorphisms of $\C$.)

1. We construct an injection $(0,\,1)\to\G_a\,,\,x\mapsto h$ as follows:
expand $x$ in binary digits (not ending on $\overline1$). Set $h:=h_n$ or
$h:=\mathrm{id}$ on $\E_n\,$, if the $n$-th digit is $1$ or $0$,
respectively, and $h:=\mathrm{id}$ on $\M\setminus\bigcup\E_n\,$. Although
the sequence of sets $\E_n$ will accumulate somewhere, continuity of $h$ can
be shown by employing $\diam(\E_n)\to0$. --- Conversely, to obtain an
injection $\G_a\to(0,\,1),\,h\mapsto x$, enumerate the hyperbolic components
$(\Omega_n)_{n\in\N}$, and denote the $n$-th prime number by $p_n\,$. Now
$x$ shall have the digit $1$ at the place $p_n^m\,$, iff
$h:\Omega_n\to\Omega_m\,$. The mapping $h\mapsto x$ is injective, since the
homomorphism from $\G_a$ to the permutation group of hyperbolic components
is injective: if $h$ is mapping every hyperbolic component to itself, it is
fixing the points of intersection of closures of hyperbolic components,
i.e.~all roots of satellite components. These are dense in $\partial\M$,
thus $h=\mathrm{id}$. --- By the two injections, $|\G_a|=|(0,\,1)|=|\R|$.

If $h_1,\,h_2\in\G_a$ with $h_1\neq h_2\,$, there is a hyperbolic component
$\Omega$ with $h_1(\Omega)\neq h_2(\Omega)$. By Prop.~\ref{PtG},
 $\mathcal{N}:=\{h\in\G_a\,|\,h(\Omega)=h_1(\Omega)\}$ is an open
neighborhood of $h_1\,$, and
 $\G_a\setminus\mathcal{N}=\bigcup\{h\in\G_a\,|\,h(\Omega)=\Omega'\}$ is an
open neighborhood of $h_2\,$, where the union is taken over all hyperbolic
components $\Omega'\neq h_1(\Omega)$. Thus $h_1$ and $h_2$ belong to
different connected components, and $\G_a$ is totally disconnected.

2. We have $d(h_n\,,\,\mathrm{id})\le2\diam(\E_n)\to0$ as $n\to\infty$, thus
$\mathrm{id}$ is not isolated in $\G_a\,$. Since composition is continuous,
no point is isolated, and $\G_a$ is perfect.

Choose a homeomorphism $h:\E_\sM\to\E_\sM$ according to Thm.~\ref{Th}, which
is expanding at $a$ and contracting at $b$, extend it by the identity to
$h\in\G_a\,$. The iterates of $h$ satisfy $h^k(a)=a$ and $h^k(c)\to b$ for
all $c\in\E_\sM\setminus\{a\}$, thus the pointwise limit of $(h^k)_{k\in\N}$
is not continuous. The sequence does not contain a subsequence converging
uniformly, and $\G_a$ is not sequentially compact, a fortiori not compact.
--- If $\mathcal{N}$ is a neighborhood of $\mathrm{id}$ in $\G_a\,$, fix an
$n$ such that $\mathcal{N}$ contains the ball of radius $2\diam(\E_n)$
around $\mathrm{id}$, then $\mathcal{N}$ contains the sequence
$(h_n^k)_{k\in\N}$. Thus $\mathcal{N}$ is not compact, and $\G_a$ is not
locally compact.

3. When $\mathcal F$ does not satisfy the finiteness condition, there is a
sequence $(h_n)\subset\mathcal F$ and a hyperbolic component $\Omega$, such
that the period of $h_n(\Omega)$ (or $h_n^{-1}(\Omega)$) diverges. Assume
$h_n\to h$, then $h_n(\Omega)=h(\Omega)$ for $n\ge N_0$ according to
Prop.~\ref{PtG}, a contradiction. If $\mathcal{F}$ satisfies the finiteness
condition, a diagonal procedure yields a subsequence which is eventually
constant on every hyperbolic component, thus respecting the partial order of
hyperbolic components. Assuming local connectivity, all fibers are trivial
\cite{sf2}, and $\lim h_n$ is obtained analogously to \cite[Sect.~9.3]{wjt}.
\mybox

Two rational angles with odd denominators are \emph{Lavaurs-equivalent}, if
the corresponding parameter rays are landing at the same root. Denote the
closure of this equivalence relation on $S^1$ by $\sim$. The \emph{abstract
Mandelbrot set} is the quotient space $S^1/\!\sim$ \cite{kkif}, it is a
combinatorial model for $\partial\M$, which will be homeomorphic to
$\partial\M$ if $\M$ is locally connected. (It is analogous to Douady's
\emph{pinched disk model} of $\M$.) Orientation-preserving homeomorphisms of
the abstract Mandelbrot set are described by orientation-preserving
homeomorphisms $\bH:S^1\to S^1$ that are compatible with $\sim$. According
to Sect.~\ref{41}, every homeomorphism $h:\E_\sM\to\E_\sM$ constructed by
surgery defines such a circle homeomorphism (extended by the identity), and
the homeomorphism group of $S^1/\!\sim$ has the properties given in
Thm.~\ref{TG}. In fact, these homeomorphisms of the abstract
Mandelbrot set can be constructed in a purely combinatorial way, without
using quasi-conformal surgery \cite{wjt}.

\small
Present address: Wolf Jung\rule{0mm}{10mm}, Inst.~Reine~Angew.~Math.\\
RWTH Aachen, D-52056 Aachen, Germany\\
\qquad\href{http://www.iram.rwth-aachen.de/~jung}%
{http://www.iram.rwth-aachen.de/$\sim$jung}\\
\href{mailto:jung@iram.rwth-aachen.de}{jung@iram.rwth-aachen.de}

\begin{thebibliography}{99} \small
\bibitem{bd} B.~Branner, A.~Douady, Surgery on complex polynomials, in:
   \emph{Holomorphic dynamics}, Proc.~2nd Int.~Colloq. Dyn.~Syst.,
   Mexico City, LNM~{\bf 1345}, 11--72 (1988).
\bibitem{bfl} B.~Branner, N.~Fagella, Homeomorphisms between limbs
   of the Mandelbrot set, J.~Geom.~Anal.~{\bf 9}, 327--390 (1999).
\bibitem{bfe} B.~Branner, N.~Fagella, Extensions of Homeomorphisms
   between Limbs of the Mandelbrot Set,
   Conform.~Geom.~Dyn.~{\bf 5}, 100--139 (2001).
\bibitem{bhjps} X.~Buff, C.~Henriksen, Julia sets in parameter spaces,
   Comm.\ Math.\ Phys.~{\bf 220}, 333--375 (2001).
\bibitem{cg} L.~Carleson, T.~W.~Gamelin, \emph{Complex dynamics},
   Springer, New York, 1993.
\bibitem{dhp} A.~Douady, J.~H.~Hubbard, On the dynamics of polynomial-like
   mappings, Ann.~Sci. \'Ecole Norm.~Sup.~{\bf 18}, 287--343 (1985).
\bibitem{TLHa} P.~Ha\"{\inodot}ssinsky, Modulation dans l'ensemble de
   Mandelbrot, in: \emph{The Mandelbrot Set, Theme and Variations},
   Tan~L.~ed., LMS Lecture Notes~{\bf 274}, Cambridge University Press 2000.
\bibitem{wjt} W.~Jung, \emph{Homeomorphisms on Edges of the Mandelbrot
   Set}, Ph.D.~thesis RWTH Aachen 2002.
\bibitem{kkif} K.~Keller, \emph{Invariant Factors, Julia Equivalences and
   the \(Abstract\) Mandelbrot Set}, LNM~{\bf 1732}, Springer 2000.
\bibitem{lvqk} O.~Lehto, K.~I.~Virtanen, \emph{Quasiconformal Mappings in
   the Plane}, Springer 1973.
\bibitem{mmu} C.~T.~McMullen, The Mandelbrot set is universal,
   in: \emph{The Mandelbrot Set, Theme and Variations},
   Tan~L.~ed., LMS Lecture Notes~{\bf 274}, Cambridge University Press 2000.
\bibitem{rscr} J.~Riedl, D.~Schleicher, On Crossed Renormalization of
   Quadratic Polynomials, in: \emph{Proceedings of the $1997$ conference
   on holomorphic dynamics}, RIMS Kokyuroku {\bf 1042}, 11--31, Kyoto 1998.
\bibitem{rt} J.~Riedl, \emph{Arcs in Multibrot Sets, Locally
   Connected Julia Sets and Their Construction by Quasiconformal
   Surgery}, Ph.D.\ thesis TU~M\"unchen 2000.
\bibitem{prt} P.~Roesch, \emph{Topologie locale des method\'es de Newton
   cubiques}, Ph.D.~thesis ENS de Lyon 1997.
\bibitem{sf2} D.~Schleicher, On Fibers and Local Connectivity of
   Mandelbrot and Multibrot Sets, IMS-preprint 98-13a (1998).
\bibitem{ser} D.~Schleicher, Rational Parameter Rays of the
   Mandelbrot Set, Ast\'erisque {\bf 261}, 405--443 (2000).
\bibitem{tls} Tan~L., Similarity between the Mandelbrot set and Julia sets,
   Comm.\ Math.\ Phys.\ {\bf 134}, 587--617 (1990).
\bibitem{tlbcn} Tan~L., Branched coverings and cubic Newton maps,
   Fundam.~Math.~{\bf 154}, 207--260 (1997).
\end{thebibliography}
\end{document}